\documentclass[11pt]{article}
\usepackage{rotating}
\usepackage{amssymb}
\usepackage{color,todonotes}
\usepackage{amsthm}
\usepackage{amsfonts}
\usepackage{graphicx}
\usepackage{graphics}
\usepackage{setspace}
\usepackage{float}
\usepackage{bm}
\usepackage{amsmath}
  {
      \theoremstyle{plain}
      
  }
\usepackage{algorithm}
\usepackage{algpseudocode}
\usepackage{appendix}
\usepackage[format=plain,justification=raggedright,singlelinecheck=false]{caption}
\usepackage[justification=centering]{caption}
\usepackage{url}
\usepackage{listings}
\usepackage{color} 
\definecolor{mygreen}{RGB}{28,10,0} 
\definecolor{mylilas}{RGB}{170,55,241}
\makeatletter
\def\BState{\State\hskip-\ALG@thistlm}
\makeatother \numberwithin{equation}{section}
\newtheorem{theorem}{Theorem}[section]

\newtheorem{proposition}[theorem]{Proposition}
\newtheorem{definition}[theorem]{Definition}
\newtheorem{corollary}{Corollary}

\newtheorem{lemma}[theorem]{Lemma}

\newtheorem{Remark}[theorem]{Remark}

\usepackage{amsthm}
\let\oldproofname=\proofname
\renewcommand{\proofname}{\rm\bf{\oldproofname}}

\newcommand{\x}{\bm{x}}

\begin{document}
\date{\small\textsl{\today}}
\vspace{-2mm}

\title{Generalized MLS and LS-SVR for Allen–Cahn Equations on Spheres}

\author{
\large Vahid Mohammadi$^{\mbox{\small \em
a}}$, \large
Mehdi Dehghan$^{\mbox{\small
\em b}}$, \large Stefano De~Marchi$^{\mbox{\small \em
	c}}$\footnote{Corresponding author. {\em E-mail addresses:}
vahid.mohammadi@sru.ac.ir (V. Mohammadi),
mdehghan@aut.ac.ir, mdehghan.aut@gmail.com (M. Dehghan), 
stefano.demarchi@unipd.it (S. De Marchi).}\\ \\
\small{\em $^{\mbox{\footnotesize a}}$Department of Mathematics, Faculty of Science,}\\
\small{\em Shahid Rajaee Teacher Training University, Tehran,
    16785-163, Iran}\\ \\
\small{\em $^{\mbox{\footnotesize b}}$Department of Applied Mathematics, Faculty of Mathematics and Computer Sciences,}\\
\small{\em Amirkabir University of Technology (Tehran Polytechnic), No. 424, Hafez Ave.,15914, Tehran, Iran}\\ \\
\small{\em $^{\mbox{\footnotesize c}}$Department of Medicine and Padova Neuroscience Center, University of Padua, Italy}\vspace{-4mm}\\
} \maketitle
\leftskip=0pt \hrule\vskip 8pt \hspace*{-.8cm} {\bfseries Abstract.}
This work extends least-squares support vector regression (LS-SVR) based on the generalized moving least squares (GMLS) approximation to the numerical solution of partial differential equations on spheres. We show that the proposed approach constructs a kernel-based approximant in its dual form and establishes its connection with the associated reproducing kernel Hilbert space (RKHS). Furthermore, solving the PDE through the proposed LS-SVR framework is shown to be equivalent to solving a smoothing problem or a kernel ridge regression problem in the corresponding RKHS.

As applications, the method is employed to solve the Allen–Cahn equation and its conservative counterpart on spheres. An implicit first-order operator-splitting scheme is used to separate each model into linear and nonlinear subproblems. The linear subproblem is discretized by the backward differentiation formula and approximated at each time step using LS-SVR derived from a GMLS-based collocation formulation, while the nonlinear part is treated analytically. The resulting scheme is shown to admit a unique numerical solution at every time step. For the conservative Allen–Cahn model, a correction procedure is incorporated to preserve mass. Numerical experiments demonstrate the accuracy and effectiveness of the proposed framework.
\bigskip
\leftskip=0pt \hrule\vskip 8pt \hspace*{-.8cm}
 \hspace*{.4cm}
\hspace*{.5cm}\\
\begin{small}$\textbf{\emph{Keywords}}$: Generalized moving least
	squares approximation; Least-squares support vector regression;
	Reproducing Kernel Hilbert space and its smoothing form; Kernel ridge regression; Allen-Cahn equation; Conservative Allen-Cahn equation; PDEs on spheres,
	Time-splitting technique.
\end{small}\\[0.2cm]
\textbf{{Mathematics Subject Classification (2020)}}:  62J02, 33C55,
65Nxx, 41Axx.
\section{Introduction}\label{sect-intro}
In this paper, we conduct a numerical investigation of two partial
differential equations (PDEs) defined on the sphere, i.e., the
Allen-Cahn (AC) equation and the conservative Allen-Cahn (CAC)
equation, using the new approach.

The AC equation on sphere $\mathbb{S}^2$ reads as follows
\cite{Allen,Choi}:
\begin{equation}\label{e1}
\dfrac{{\partial u (\bm{x},t)}}{{\partial t}} ={\Delta
_{\mathbb{S}^2}}u (\bm{x},t)  - \dfrac{{F^{\prime}(u
(\bm{x},t))}}{{{{\epsilon} ^2}}},
\end{equation}
in which $\bm{x} \in \mathbb{S}^2$, $t \in (0,T]$ and
 $\Delta_{\mathbb{S}^2}$ shows the Laplace-Beltrami
operator, and $u (\bm{x},t)$ is the phase-field variable (order parameter). Also,
$F(u(\bm{x},t))=0.25(u(\bm{x},t)^{2}-1)^{2}$ shows the free-energy function.
$\epsilon$ shows the gradient energy coefficient corresponding to
the interfacial energy.

Besides, the CAC equation on sphere $\mathbb{S}^2$ is formulated by
\cite{kim-CAC}:
\begin{equation}\label{e1-CAC}
    \dfrac{{\partial u (\bm{x},t)}}{{\partial t}} ={\Delta _{\mathbb{S}^2}}u
(\bm{x},t)  - \dfrac{{F^{\prime}(u (\bm{x},t))}}{{{{\epsilon}
^2}}}+\beta(t)\sqrt{F(u(\x,t))},
\end{equation}
in which $\sqrt{F(u(\bm{x},t))}=0.5|u(\bm{x},t)^2-1|$, and
\begin{equation}\label{beta-def}
\beta(t)=\dfrac{\displaystyle \int_{\mathbb{S}^2}F^{\prime}(u(\x,t))
~d\bm{x}}{\epsilon^2 \displaystyle \int_{\mathbb{S}^2}\sqrt{F(u(\x,t))}~d\bm{x}}.
\end{equation}
The Lagrange multiplier $\beta(t)$ is used to satisfy the mass conservation
property for variable $u(\bm{x},t)=u$ in time \cite{Kim11,kim-CAC}. 
As a note, the denominator of $\beta(t)$ is  positive.  Indeed,  $\sqrt{F(u)}=0.5|u^2-1|$ is nonnegative for $u \in [-1,1 ]$, and it vanishes only at phases $u=-1$ or $u=1$.  In the proposed diffuse-interface model, the transition layer between phases has finite thickness $\epsilon$. Hence, the region where $-1<u<1$ has positive measure,  making the integral strictly positive. The only case where this term could be zero occurs at the trivial homogeneous state, i.e., $u=-1$ or $u=1$ everywhere, where the diffusion term is absent.  
Using an asymptotic analysis by considering the diffusion term  gives a tanh profile as the solution showing  a diffuse-interface of  positive width always exists (cf. \cite{Mohammadi-AMM-2022}). 

In Appendix A, we show that the phase-field variable $u$ of model (\ref{e1-CAC}) preserves mass in time.

Until now, numerous numerical methods have been developed for
solving PDEs on surfaces (to predict natural phenomena,
such as image inpainting \cite{Feng,Li-surface}, pattern formations
\cite{Fuselier}, and tumor growth \cite{Wu}) such as \eqref{e1} and
\eqref{e1-CAC}, as the availability of exact solutions is impossible in most cases. In previous research works related to solving Eqs.
\eqref{e1} and \eqref{e1-CAC} several numerical approaches were
developed, including the classical finite difference (FD) method on
the narrow band domain boundary, combined with an explicit
time-splitting scheme of order one for solving the AC equation on
surfaces \cite{Choi}. Moreover, a CAC equation on surfaces was
numerically investigated by utilizing the standard FD scheme in
space and a first-order time-splitting idea in time \cite{Kim11}.
Authors of \cite{MohDehMir} used an explicit time-splitting method
of order one, combined with the restricted RBFs interpolant to solve
the AC equation on $2$-dimensional surfaces. Recently, a new numerical solution to the CAC equation on the cubic surface has been
found in \cite{Hwang}.

Overall, the numerical techniques applied to the surface PDE
problems fall into two groups, i.e., the intrinsic and embedded
narrow band methods (cf.
\cite{Choi,Flyer-Wright:2007-1,Flyer-Wright:2009-1,Fuselier,Gera,Kim11,
Kimsurface,Lehto1,mirzaei:2017-2,MohDehMir, morton-neamtu:2002-1,
Shankar,xiao-et-al:2017}). The first group of these methods attempts
to find an approximation to the differential operator in the
original dimension. Of course, this approach faces difficulty with
the inherent coordinate distortions or singularities arising in the
metric terms of the differential operators (cf. \cite{Fuselier}). In
the second group, the differential operators are extended to the
extrinsic coordinates, which avoids all coordinate singularities.
Accordingly, any standard numerical methods such as the finite
difference (FD) schemes on $3D$ regular nodes and the finite element
method (FEM) on unstructured meshes can be implemented. Of course,
these methods increase the computational costs because they solve
the problem in a dimension at least one greater than the original
dimension of the main equation. Additionally, it needs to impose
artificial boundary conditions, which may lead to the low accuracy
of the applied numerical method. Hence, to avoid all the
disadvantages of the mentioned approaches, a high-order numerical
method applied to PDE problems on surfaces has been introduced in
\cite{Fuselier} that is based on the radial basis functions (RBFs)
interpolation, and it integrates the benefits of both intrinsic and
embedded narrow-band techniques (see also \cite{MohDehMir}).
Besides, its recent developments, i.e., the RBF partition of unity
(RBF-PU) method \cite{Mir-D-RBF-PU} and the RBF finite difference
(RBF-FD) scheme \cite{Bayona-2024, Shankar} have successfully been
applied to the surface PDEs. In addition, the spherical basis
functions (SBFs) have been applied to find the numerical solution of
the PDEs on spheres (cf.
\cite{Gia-2004,Gia-2005,Gia-2012,mirzaei:2017-2,Narcowich-2017,Xu}.

Another possible numerical method on the spheres is to use
 zonal basis functions, which were investigated in \cite{Hubbert, Narcowich-2002}.
 A zonal kernel is constructed by the restricted original positive definite radial kernel, say $\Phi=\phi(||.||_{2}):\mathbb{R}^{d+1} \to \mathbb{R}$ to $\mathbb{S}^d=\{\bm{x} \in \mathbb{R}^{d+1}:\, ||\bm{x}||_{2}=1\}$ ($d$-dimensional sphere) such that for any $\bm{x}, \bm{y} \in \mathbb{S}^d$, we have $||\bm{x}-\bm{y}||_{2}=\sqrt{2-2\bm{x}^{T}\bm{y}}$, and thus $\Phi(\bm{x},\bm{y})=\phi(\sqrt{2-2\bm{x}^{T}\bm{y}}):=\psi(\bm{x}^{T}\bm{y})$.  
 Also, these basis functions have been used in the Petrov-Galerkin form of
  the diffusion equations on the unit sphere as the space of trial functions
  (cf. \cite{mirzaei:2017-2,AhmadiDarani-Mirzaei}). In addition, in \cite{Jones-2023}, by applying the
local coordinates of the tangent plane, the differential operators
on the two-dimensional surfaces have been approximated by the RBF-FD
and the GMLS procedures.

On the other hand, one of the famous meshfree approaches is the
generalized moving least squares (GMLS) approximation on spheres
introduced in \cite{mirzaei:2017-1}, known as a direct
approximation. Comparing this method with the original moving least
squares (MLS) technique presented in \cite{Wendland-MLS-Sphere}, the GMLS
technique does not need to compute the derivatives (in many PDE
problems on spheres, the surface gradient, and the Laplace-Beltrami
operators should be computed) of the complicated shape functions. As
will be observed, the approximate solution is found in the space
spanned by the spherical harmonics of degree less than or equal to
$m$ (the maximum degree of spherical harmonics used here), and the Laplace-Beltrami operator is approximated easily due
to the main property of the spherical harmonics (see Section
\ref{GMLS-spheres}). Besides, the proposed method does not depend on
any mesh or background cell for the approximation, similar to the
meshfree methods.

In this paper, we propose a new numerical procedure based on
a supervised learning algorithm, i.e., the least-squares support
vector regression (LS-SVR)
, considering the GMLS approximation \cite{mirzaei:2017-1} to solve both AC and CAC equations on
 $\mathbb{S}^2$. Moreover, we show that the constructed kernel function on the sphere induces a reproducing kernel Hilbert space (RKHS) in the dual representation of the approximate solution proposed in LS-SVR
  Furthermore, we prove that the original minimization problem in the LS-SVR is equivalent to a smoothing problem or kernel ridge regression (KRR) in its corresponding RKHS. As far as we are aware, in the literature, this point of view related to the proposed supervised learning method is not considered at all. Hence, this is one of the main reasons for writing this research work. In what follows, we only mention some recent research works on this topic, where the proposed supervised learning algorithm has been applied to solve partial (ordinary) differential equations or integral equations.

The LS-SVR is a type of the least-squares support vector machine
(LS-SVM), which is applied for function estimation
\cite{Mehrkanoon-2012, Razzaghi-SVM}. This approach is known as one
of the main categories of
 supervised learning \cite{Torkian}. In \cite{Mehrkanoon-2014, Razzaghi-SVM}, this idea is successfully
 applied to solve both ordinary and partial differential equations. Also, the authors of \cite{Parand-1} have
 extended the LS-SVR to solve the Volterra integral equations. In \cite{Parand-2}, this procedure has been
  developed to find the numerical solutions to Fredholm integral equations. The same approach, combined with
   spectral methods, is employed for solving numerically differential equations on unbounded domains \cite{Pakniyat}.
    In \cite{Abbaszadeh-JCP}, a new algorithm based on a combination of the reduced-order LS-SVR with the isogeometric
     collocation method has been applied to solve the Cahn-Hilliard-Navier-Stokes (CHNS) equations. See also
      \cite{Abbaszadeh-MCS}, where the LS-SVR combined with a meshfree approximation has been developed to
      find the numerical solution of the fractional PDEs on diverse complex geometries.

Main contributions of this paper are itemized as follows:
\begin{itemize}

  \item
  We show that the proposed supervised learning algorithm on spheres can be interpreted as a smoothing problem or a KRR in its corresponding RKHS. 

    \item
    For the first time, a supervised learning algorithm, i.e., the LS-SVR based on the GMLS approximation is developed to solve numerically two nonlinear phase-field problems on spheres.
    The applied GMLS technique allows us to approximate the Laplace-Beltrami operator more easily than
    a mesh-dependent method, for example, methods developed in \cite{Choi, Kim11}.

    \item
An implicit time-splitting approach combined with the proposed learning algorithm makes a simple 
fully discrete scheme for each phase-field model. Indeed, instead of solving a nonlinear system per time step, an implicit linear and a nonlinear problems are solved, where the nonlinear part is treated analytically. As a result, no nonlinear
iterative solvers, such as the Newton-Raphson algorithm, are
required to obtain the approximate solution for the fully implicit
schemes of the studied mathematical models.

   \item
 By using the proposed supervised learning algorithm, the fully discrete problem related to the linear part of each equation is a minimization problem with constraints obtained from the collocation technique of the GMLS method at each time step. Accordingly, a dual form of this regression problem is equivalent to solving a linear system of algebraic equations with a positive definite coefficient matrix, resulting in a unique approximate intermediate solution.
\end{itemize}
Additionally, unlike the methods presented in \cite{Choi, Kim11}, the developed numerical technique does not require defining any quasi-Neumann boundary conditions on the narrow band domain boundary. On the other hand, the presented LS-SVR is meshfree because no triangulation or background cell is required for constructing the trial approximation via the GMLS procedure. Instead, it only depends on the number of nodes and their locations scattered on spheres. This makes it easy to execute on the PDE problems.

The other parts of the paper are organized as follows. In Section
\ref{GMLS-spheres}, the GMLS technique on spheres is explained
briefly. Section \ref{sec-LS-SVR} proposes the implementation of the
LS-SVR to solve the PDE problem on spheres. In the next section, we first prove that the space of real-valued functions belonging to an approximation space is an RKHS. Then,
we show how the proposed LS-SVR can be interpreted in the sense of a smoothing problem or a KRR in its associated RKHS.
In Section \ref{sect-LS-SVR-AC-CAC}, we first apply an
implicit time-splitting scheme of order one to Eqs. (\ref{e1}) and
(\ref{e1-CAC}) for time discretization. Then, we develop the
proposed LS-SVR based on the GMLS approximation to derive the fully
discrete schemes for both phase-field models.
In the next section,
 the ability of the presented numerical procedure is demonstrated through the investigation of
  some simulation results.  Finally, we give a conclusion in Section
\ref{sect-Conclusion}.

\section{The GMLS procedure on spheres} \label{GMLS-spheres}
In this section, a brief description of the GMLS technique on a 2D
sphere, proposed in \cite{mirzaei:2017-1}, is presented. As is
already mentioned, the method is meshfree. This means that any
background cell or triangulation on spheres is not necessary for the
approximation. Instead, the scattered points, like
$X=\{\bm{x}_{1},...,\bm{x}_{N}\}$, on spheres are used to construct
the trial approximation. However, a sequence of such sets should be
quasi-uniform, meaning there exists a constant $c_{qu}>1$
independent of $X$ for all sets of this sequence, such that the mesh-ratio
$\rho_{X}:={h_{X,\mathbb{S}^2}}/{q_{X}} \leq c_{qu}$, where
$h_{X,\mathbb{S}^2}$ and $q_{X}$ are the fill distance and
separation distance of points $X$, respectively (cf.
\cite{mirzaei:2017-1}).

We are interested in finding an approximation to the
Laplace-Beltrami operator $\Delta_{\mathbb{S}^{2}}:=\texttt{L}$,
applied to the function $u \in C^{m+1}(\mathbb{S}^2)$ for some $m
\geq 0$. To this end, the GMLS approximation to
$\texttt{L}u(\bm{x})$ based on the points $X$ represented by
$s^{\texttt{L}}_{X}(\bm{x})$, where $\bm{x} \in \mathbb{S}^2$, is
given by \cite{mirzaei:2017-1}
\begin{equation}\label{GMLS-Sphere1}
s^{\texttt{L}}_{X}(\bm{x})=\displaystyle \sum_{j\in
I(\bm{x})} a^{\texttt{L}}_{j}(\bm{x})u_{j}=\bm{a}^{\texttt{L}}(\bm{x})^T \bm{u}_{X},
\end{equation}
in which the coefficients $\{a^{\texttt{L}}_{j}(\bm{x}): j \in I(\bm{x})\}$ are computed
through solving the weighted minimization problem over the spherical harmonics space with
degree $m$ or less, say $\pi_{m}(\mathbb{S}^2)$ with dimension $N(3,m)$
(cf. \cite[Theorem 4.1]{mirzaei:2017-1}). Note that
$$N(3,0)=1,\quad N(3,m)=\frac{(2m+2)\Gamma(m+2)}{\Gamma(m+1)\Gamma(3)},\quad m \geq 1.$$
Besides, $\bm{u}_{X}=[u(\bm{x}_{i})]_{1 \leq i \leq N}$ denotes
 the value of $u(\bm{x})$ at points $X$. Besides, a set of indices, i.e., $I(\bm{x})$ is defined by
$$ \{1 \leq j \leq N: dist(\bm{x},\bm{x}_{j})\leq \delta\},$$
 where $dist(\bm{x},\bm{y})=arccos(\bm{x}^T\bm{y})$ denotes the geodesic
distance between two points $\bm{x}, \bm{y}$ on $\mathbb{S}^2$, and $0<\delta < \pi$ is the radius of the spherical
 cap \cite{mirzaei:2017-1,Wendland-MLS-Sphere}. By assuming the unisolvency of $X$ on each spherical cap, we have
\begin{equation}\label{GMLS-Sphere6}
\bm{a}^{\texttt{L}}(\bm{x})=W(\bm{x})P\Big(P^{T}W(\bm{x})P\Big)^{-1}\texttt{L}Y(\bm{x}),\quad
\bm{x} \in \mathbb{S}^{2},
\end{equation}
where $P$ is a ${|I(\bm{x})| \times N(3,m)}$ matrix with
elements $Y_{\ell k}(\bm{x}_{j})$. $W(\bm{x})$ is a diagonal
matrix of size $|I(\bm{x})|$ carrying the weights
$w(\bm{x},\bm{y}):=\phi
\big(\frac{dist(\bm{x},\bm{y})}{\delta}\big)$ on its diagonal,
where $w: \mathbb{S}^{2}\times \mathbb{S}^{2} \rightarrow
[0,\infty)$ is a continuous function. This function is typically
compactly supported radial, i.e., $\phi:[0,\pi) \rightarrow
[0,\infty)$ such that for $0 \leq r \leq 1$, $\phi(r)>0$ and for
$r>1$, $\phi(r)=0$ (of course, it
is possible to utilize the weight function such that when $r\to 0$,
$\phi \to \infty$) \cite{mirzaei:2017-1,Wendland-MLS-Sphere}.
Note that, in (\ref{GMLS-Sphere6}), the matrix $A(\bm{x}):=P^{T}W(\bm{x})P$ is invertible because it is positive definite (cf. \cite{mirzaei:2017-1}).

Besides, in relation \eqref{GMLS-Sphere6}, $[\texttt{L}Y_{\ell^{'}k^{'}}(\bm{x})]_{0 \leq  \ell^{'} \leq m,\,1 \leq k^{'} \leq N(2,\ell^{'})}={\texttt{L}}Y(\bm{x}) \in \mathbb{R}^{N(3,m)}$ is the vector of the spherical
harmonic functions with degree less than or equal $m$, where $\texttt{L}$ acts on the vector
$Y(\bm{x})$ that can be computed throughout the following relation (cf. \cite{Atkinson-2,Dai,Dai2}).
\begin{equation}\label{GMLS-Sphere2}
\texttt{L}Y(\bm{x})=-\ell(\ell+1)Y(\bm{x}).
\end{equation}
Finally, an approximation to $\texttt{L}u(\bm{x})$ at points $X$ via the GMLS approximation is
\begin{equation}\label{GMLS-Sphere-new-2}
\texttt{L}u(\bm{x})|_{X}\approx\bm{a}^{\texttt{L}^{T}}(\bm{x})|_{X} \bm{u}_{X}= L_{X}\bm{u}_{X},
\end{equation}
in which $L_{X}=[a^{\texttt{L}}_{j}(\bm{x}_{i})]_{1\leq i,j \leq N}$ is a sparse differentiation matrix, and its elements are computed according to relation (\ref{GMLS-Sphere6}).

On the other hand, if we consider the
point-evaluation functional, $\delta_{\bm{x}}$ instead of
 $\texttt{L}$, then the GMLS approximation to the function $u$ represented by $s_{X}(\bm{x})$,
 where $\bm{x} \in \mathbb{S}^2$, is given by
\begin{equation}\label{GMLS-Sphere-new}
s_{X}(\bm{x})=\displaystyle \sum_{j\in
I(\bm{x})} a_{j}(\bm{x})u_{j}=\bm{a}^{T}(\bm{x}) \bm{u}_{X},
\end{equation}
which the coefficients $a_{j}(\bm{x})$ satisfy the delta-kronecker property,
i.e., if $i=j$, $a_{j}(\bm{x}_{i})=1$, else-if $i\neq j$, $a_{j}(\bm{x}_{i})=0$, where $i,j=1,...,N$.
This means that the approximation \eqref{GMLS-Sphere-new} is an interpolant form.
Also,
\begin{equation}\label{GMLS-Sphere-new-1}
\bm{a}(\bm{x})=W(\bm{x})P\Big(P^{T}W(\bm{x})P\Big)^{-1}Y(\bm{x}),\quad \bm{x} \in \mathbb{S}^2.
\end{equation}
 Also, $u(\bm{x})$
can be approximated at points $X$ using the proposed GMLS idea as follows:
\begin{equation}\label{GMLS-Sphere-new-2-2}
u(\bm{x})|_{X}\approx\bm{a}^T(\bm{x})|_{X} \bm{u}_{X}= A_{X}\bm{u}_{X},
\end{equation}
where $A_{X}=[a_{j}(\bm{x}_{i})]_{1\leq i,j \leq N}$ is a sparse matrix. 

As a note, this idea is referred to as the direct approximation due to acting the linear
functional $\texttt{L}$ (or any other linear functionals) only on spherical harmonics,
as claimed in \cite{mirzaei:2017-1} and seen here. The interested reader can also refer
to \cite{mirzaei:2017-1} for finding most materials corresponding to the theoretical aspects of this approximation.

\section{The LS-SVR on spheres}\label{sec-LS-SVR}
In this section, we extend the LS-SVR \cite{Mehrkanoon-2012} to
spheres for solving the PDE problems (or in the special case, the approximation problem) in the sense of the collocation
form. 

To start, assume that the linear differential equation is given by
\begin{align}\label{LS-SVR-1}
    \mathcal{L}u(\bm{x})=f(\bm{x}),\quad \bm{x} \in \mathbb{S}^2,
\end{align}
in which $\mathcal{L}$ denotes a differential operator defined on
$\mathbb{S}^2$, and $u$ and $f$ are sufficiently smooth functions. Now,
consider the training (trial) data $\{\bm{x}_1,...,\bm{x}_{N}\}$ on
$\mathbb{S}^2$. Also, the solution to the regression problem is
considered as follows:
\begin{align}\label{approximate-LS-SVR}
\hat{u}(\bm{x})=\bm{a}^{T}(\bm{x})\bm{w}=\displaystyle \sum_{j=1}^{N}a_{j}(\bm{x})w_{j},
\end{align}
in which $\bm{w}=[w_{1},...,w_{N}]^{T}$ and
$\bm{a}^{T}(\bm{x})=[a_{1}(\bm{x}),...,a_{N}(\bm{x})]$ contain the
weights and the basis functions, respectively.
To find the weights, we should solve the following minimization problem (cf.
\cite{Mehrkanoon-2012, Razzaghi-SVM}).
\begin{align}
    &\min_{\bm{w},\bm{e} \in \mathbb{R}^{N}} \dfrac{1}{2}\bm{w}^{T}\bm{w}+\dfrac{\gamma}{2}\bm{e}^{T}\bm{e},\label{LS-SVM-1}\\
    &\mathcal{L}\hat{u}(\bm{x}_{i})-f(\bm{x}_{i})=e_{i},\label{LS-SVM-2}\quad i=1,...,N,
\end{align}
where $\bm{e}=[e_{1},...,e_{N}]^{T}$ is the vector containing error value at each point $\bm{x}_{i}$
due to the proposed constraints in (\ref{LS-SVM-2}). Additionally, $\gamma$ represents the positive constant,
which is usually chosen arbitrarily by the user.
However, in the next section, we show how this parameter can be chosen based on the theory of this approach. 
 According to \cite[Theorem $7.1$]{Razzaghi-SVM}
(see also \cite{Mehrkanoon-2012}), where the Karush-Kuhn-Tucker (KKT) optimality conditions are applied,
we find the solution of minimization problem \eqref{LS-SVM-1}-\eqref{LS-SVM-1} as follows:
\begin{align}\label{LS-SVM-1-1}
    \bm{w}=A^{T}\bm{\alpha},\quad \bm{e}=-\dfrac{1}{\gamma}\bm{\alpha},
\end{align}
where $A=:A_{L}=[\mathcal{L}a_{j}(\bm{x}_{i})]_{1\le i,j\le N}$ thanks to
the GMLS approximation given in the previous section. Also, $\bm{\alpha}=:\bm{\alpha}(\gamma)$ (i.e., the vector $\bm{\alpha}$ depends on parameter $\gamma$)
is the vector of Lagrangian coefficients in the dual form of the
minimization problem (\ref{LS-SVM-1})-(\ref{LS-SVM-2}), which is
computed by solving the following linear algebraic equations
\cite{Razzaghi-SVM}.
\begin{align}\label{LS-SVM-3}
    (AA^{T}+\frac{1}{\gamma}I_{N})\bm{\alpha}=\bm{b},
\end{align}
in which $\bm{b}=[f(\bm{x}_{i})]_{1\leq i \leq N}$, and it is known as a dual form of problem
 (\ref{LS-SVM-1})-(\ref{LS-SVM-2}). Note also that, the matrix $B=:AA^T+\frac{1}{\gamma}I_{N}$ is clearly symmetric, i.e.,
 $B^{T}=B$.
 \begin{theorem}\label{Th-LS-SVR-1}
    The problem (\ref{LS-SVM-1})-(\ref{LS-SVM-2}) has a unique solution (\ref{LS-SVM-1-1}).
    \begin{proof}
    It is sufficient to show that the coefficient matrix $B$ in (\ref{LS-SVM-3}) is invertible. To do so, we show the matrix $B$ is positive definite. Indeed, for any nonzero $\bm{x} \in \mathbb{R}^{N}$,
$$\bm{x}^{T}B\bm{x}=\bm{x}^{T}(AA^T+\frac{1}{\gamma}I_{N})\bm{x}=||A^{T}\bm{x}||^{2}_{2}+\frac{1}{\gamma}||\bm{x}||^{2}_{2}>0.$$
 Hence, the matrix $B$ is invertible, which gives a unique solution $\bm{\alpha}$ for any $\gamma >0$ and also problem (\ref{LS-SVM-1})-(\ref{LS-SVM-2}) according to (\ref{LS-SVM-1-1}).
    \end{proof}
  \end{theorem}
The following lemma shows the formulation of approximate solution (\ref{approximate-LS-SVR}) in the dual form.
\begin{lemma}\label{Lemma1}
    For any $\gamma >0$, the dual representation of the approximate solution (\ref{approximate-LS-SVR}) is given by:
    \begin{align}\label{approximate-LS-SVR-Dual}
    \hat{u}(\bm{x})=\displaystyle \sum_{j=1}^{N}K(\bm{x},\bm{x}_{j})\alpha_{j}=\bm{k}^{T}(\bm{x})\bm{\alpha}:=\hat{u}_{\gamma}(\bm{x}),\quad \bm{x} \in \mathbb{S}^2,
    \end{align}
    where $K(.,.)$ is the kernel function defined by $K(\bm{x},\bm{x}_{j})=:\langle\bm{a}(\bm{x}),\mathcal{L}\bm{a}(\bm{x}_{j})\rangle$ for $j=1,...,N$, and $\langle .,. \rangle$ denotes the usual inner-product between two vectors. Besides, $\bm{k}^{T}(\bm{x})=[K(\bm{x},\bm{x}_{1})\ldots K(\bm{x},\bm{x}_{N})]$ and $\bm{\alpha}=[\alpha_{1} \ldots \alpha_{N}]$.
    \begin{proof}
        According to the first relation in (\ref{LS-SVM-1-1}), the approximate solution (\ref{approximate-LS-SVR}) can be written as follows:
        \begin{align}\label{approximate-dual-LS-SVR}
        \hat{u}(\bm{x})&=\bm{a}^{T}(\bm{x})\bm{w}\nonumber\\
        &=\bm{a}^{T}(\bm{x})A^{T}\bm{\alpha}\nonumber\\
        &=\bm{a}^{T}(\bm{x}) [\mathcal{L}\bm{a}(\bm{x}_{1}) \ldots \mathcal{L}\bm{a}(\bm{x}_{N})]\bm{\alpha}\nonumber\\
        &= [\bm{a}^{T}(\bm{x})\mathcal{L}\bm{a}(\bm{x}_{1}) \ldots \bm{a}^{T}(\bm{x})\mathcal{L}\bm{a}(\bm{x}_{N})]\bm{\alpha} \nonumber \\
        &=[\langle\bm{a}(\bm{x}),\mathcal{L}\bm{a}(\bm{x}_{1})\rangle \ldots \langle\bm{a}(\bm{x}),\mathcal{L}\bm{a}(\bm{x}_{N})\rangle ]\bm{\alpha}.
        \end{align}
        Now, by defining $\langle\bm{a}(\bm{x}),\mathcal{L}\bm{a}(\bm{x}_{j})\rangle:=K(\bm{x},\bm{x}_{j}),\,j=1,...,N$, the result is concluded.
    \end{proof}
\end{lemma}

\section{An RKHS}
Suppose in (\ref{LS-SVR-1}) $\mathcal{L}=\texttt{id}$,  that is the identity operator.  So, 
the constructed kernel in the previous section allows us to define a new Hilbert space with reproducing kernel.  Accordingly,
the new function, i.e.,  
\begin{align}\label{New-Kernel}
K(\bm{x},\bm{y}):=\displaystyle \sum_{i=1}^{N}a_{i}(\bm{x})a_{i}(\bm{y})=\langle \bm{a}(\bm{x}),\bm{a}(\bm{y})\rangle_{\mathbb{R}^{N}},\quad \forall \bm{x},\bm{y} \in \mathbb{S}^{2},
\end{align}
can be constructed.

Now, we assume that for any distinct points $\bm{x}_{1},\ldots \bm{x}_{N}$ on $\mathbb{S}^{2}$, the matrix $A$ has  full rank. Accordingly, the matrix $K:=AA^T$ is positive definite. Hence, based on \cite{Aronszajn}, there exists a reproducing kernel Hilbert space (RKHS), where $K$ is reproducing kernel.

To show this fact, let $\bm{a}:\mathbb{S}^{2} \to \mathbb{R}^{N}$ be a feature map into a Hilbert space $(\mathbb{R}^{N},\langle.,.\rangle_{\mathbb{R}^{N}})$. We consider the following Hilbert space.
\begin{align}\label{Hilbert-Space-1}
H=:H(\mathbb{S}^2):=\left\{ {{f_{\bm{w}}}:\mathbb{S}^{2}  \to \mathbb{R}:\,\,\exists \bm{w} \in \mathbb{R}^{N}\,\,\textsl{with} \,\,{f_{\bm{w}}}(\bm{x}) = {{\left\langle {\bm{w},\bm{a}(\bm{x})} \right\rangle }_{\mathbb{R}^{N}}}},\,\,\forall \bm{x} \in \mathbb{S}^{2}\right\},
\end{align} 
with the inner product defined by:
\begin{align}\label{Hilbert-Space-2}
\forall f_{\bm{w}}, f_{\bm{v}} \in H, \quad \langle f_{\bm{w}},f_{\bm{v}} \rangle_{H}:=\langle \bm{w},\bm{v} \rangle_{\mathbb{R}^{N}}.
\end{align}
Besides, its associated norm is 
\begin{align}\label{Hilbert-Space-3} ||f_{\bm{w}}||_{H}=||\bm{w}||_{\mathbb{R}^{N}}.
\end{align}
Based on the above assumption and definition, 
the map $T:\mathbb{R}^{N} \to H$ with $T(\bm{w})=f_{\bm{w}}$ for all $\bm{w} \in \mathbb{R}^{N}$, where $f_{\bm{w}}(\bm{x})={{\left\langle {\bm{w},\bm{a}(\bm{x})} \right\rangle }_{\mathbb{R}^{N}}},\,\, \forall \bm{x} \in \mathbb{S}^{2}$, is an isometric isomorphism. 
This confirms  the  well-definedness of  RKHS inner product.

We are now interested to show $H$ is a reproducing kernel Hilbert space (RKHS), where its kernel is given by (\ref{New-Kernel}).
\begin{theorem}\label{Theorem-RKHS-1}
	The function (\ref{New-Kernel}) is reproducing kernel for $H$.
	\begin{proof}
	Based on \cite[Definition $13.1$]{Fasshauer}, we show $K(.,\bm{x}):\mathbb{S}^{2} \to \mathbb{R}$ belongs to $H$ for all $\bm{x} \in \mathbb{S}^{2}$, and $f_{\bm{w}}(\bm{x})=\langle f_{\bm{w}},K(.,\bm{x})\rangle_{H}$ for all $f_{\bm{w}} \in H$ and all $\bm{x} \in \mathbb{S}^{2}$.
	
	For the first claim, assume that $\bm{x} \in \mathbb{S}^{2}$ is arbitrary and fixed. For any $\bm{x}^{'} \in \mathbb{S}^{2}$, $K(\bm{x}',\bm{x})=\langle \bm{a}(\bm{x}),\bm{a}(\bm{x}^{'}) \rangle_{\mathbb{R}^{N}}$.
	On the other hand, if we consider $\bm{w}=\bm{a}(\bm{x})$, then for all $\bm{x}^{'} \in \mathbb{S}^{2}$
	\begin{align}
	f_{\bm{w}}(\bm{x}^{'})&=\langle \bm{w},\bm{a}(\bm{x}^{'})\rangle_{\mathbb{R}^{N}}\nonumber \\
	&=\langle \bm{a}(\bm{x}),\bm{a}(\bm{x}^{'})\rangle_{\mathbb{R}^{N}}\nonumber \\
	&=K(\bm{x}',\bm{x}),
	\end{align}
	which shows $k(.,\bm{x}) \in H$.
	
	For the second claim, from the above, we know that $f_{\bm{a}(\bm{x})}(.)=K(.,\bm{x})$. Hence, for any $f_{\bm{w}} \in H$
	\begin{align}
	\langle f_{\bm{w}},K(.,\bm{x})\rangle_{H}&=\langle f_{\bm{w}},f_{\bm{a}(\bm{x})}\rangle_{H}\nonumber \\
	&=\langle \bm{w},\bm{a}(\bm{x})\rangle_{\mathbb{R}^{N}}\nonumber \\
	&=f_{\bm{w}}(\bm{x}),
	\end{align}
	which completes the proof.
	\end{proof}
\end{theorem}
There are other properties related to reproducing kernel which are given in the following proposition. 
\begin{proposition}
	For a Hilbert space of functions with reproducing kernel $K$ such as $H$, we have
	\begin{itemize}
		\item $\forall \bm{x}, \bm{y} \in \mathbb{S}^{2}$, $K(\bm{x},\bm{y})=\langle K(.,\bm{x}),K(.,\bm{y})\rangle_{H}$.
		\item 
		$\forall \bm{x}, \bm{y} \in \mathbb{S}^{2}$, $K(\bm{x},\bm{y})=K(\bm{y},\bm{x})$.
		\item 
		if $\mathop {\lim }\limits_{n \to \infty } {\left\| {f - {f_n}} \right\|_H} = 0$, then $\mathop {\lim }\limits_{n \to \infty } \left| {f(\bm{x}) - {f_n}(\bm{x})} \right| = 0$, for all $\bm{x} \in \mathbb{S}^{2}$.
	\end{itemize}
\begin{proof}
	The first property can be obtained directly from Theorem \ref{Theorem-RKHS-1}. The second property is concluded from (\ref{New-Kernel}). To prove the last property, assume that $\bm{x} \in \mathbb{S}^{2}$ is arbitrary. Using the famous Cauchy-Schwarz inequality, we can write
	\begin{align}
\left| {f(\bm{x}) - {f_n}(\bm{x})} \right| &= \left| {{{\left\langle {f - {f_n},K(.,\bm{x})} \right\rangle }_H}} \right|\nonumber \\
& \le {\left\| {f - {f_n}} \right\|_H}{\left\| {K(.,\bm{x})} \right\|_H}\nonumber \\
&={\left\| {f - {f_n}} \right\|_H}{\left\| {\bm{a}(\bm{x})} \right\|_{{\mathbb{R}^{N}}}}.
	\end{align}
	Accordingly, if $n \to \infty$, then the desired result can be obtained, and it completes the proof. 
\end{proof}
\end{proposition}
We now consider the following pre-Hilbert space, namely
 $H_{pre}$.
\begin{align}\label{Pre-Hilbert space-1}
H_{pre}:=\{\displaystyle \sum_{i=1}^{N}\alpha_{i}K(.,\bm{x}_{i}):\, N \in \mathbb{N},\, \alpha_{1},...,\alpha_{N} \in \mathbb{R},\,\bm{x}_{1},...,\bm{x}_{N} \in \mathbb{S}^2\}.
\end{align}

\begin{proposition}\label{Pre-Hilbert space-3}
	Let $\bm{a}:\mathbb{S}^2 \to \mathbb{R}^{N}$ be a feature map into a Hilbert space $(\mathbb{R}^{N},\langle.,.\rangle_{\mathbb{R}^{N}})$. Then, we have 
	\begin{itemize}
		\item The pre-Hilbert space $H_{pre}$ is dense in $H$. 
		\item 
$\forall f,g \in H_{pre},\quad \langle f,g \rangle_{H}=\displaystyle \sum_{j=1}^{N}\sum_{m=1}^{N} \alpha_{j} \beta_{m}k(\bm{x}_{j},\bm{x}_{m})=\bm{\alpha}^{T}K\bm{\beta}.$
       \item
	$\forall f \in H_{pre},\quad ||f||^{2}_{H}=\displaystyle \sum_{j=1}^{N}\sum_{m=1}^{N} \alpha_{j} \alpha_{m}k(\bm{x}_{j},\bm{x}_{m})=\bm{\alpha}^{T}K\bm{\alpha}.$
	\end{itemize}
	\begin{proof}
		The first part can be proved in the same manner presented in \cite[Theorem $4.21$]{Steinwart}. For the second part, for any $ f,g \in H_{pre}$, we can write 
		$$f(.)=\displaystyle \sum_{j=1}^{N} \alpha_{j}K(.,\bm{x}_{j}),\quad g(.)=\displaystyle \sum_{m=1}^{N} \beta_{m}K(.,\bm{x}_{m}).$$
		From the reproducing property of $H$, we can write
		\begin{align}
		 \langle  f,g\rangle_{H}
		 &=\langle \displaystyle \sum_{j=1}^{N} \alpha_{j}K(.,\bm{x}_{j}), \displaystyle \sum_{m=1}^{N} \beta_{m}K(.,\bm{x}_{m})\rangle_{H} \nonumber \\
		 &=\displaystyle \sum_{j=1}^{N} \displaystyle \sum_{m=1}^{N} \alpha_{j}\beta_{m} K(\bm{x}_{j},\bm{x}_{m}) \nonumber \\
         &=\bm{\alpha}^{T}K\bm{\beta}.
		 \end{align}
		 Finally, the last part is obtained from the second part, where $g$ is replaced by $f$.
	\end{proof}
\end{proposition}
\begin{Remark}
   It should be noted that if we want to have an infinite-dimensional space, the Hilbert space $\ell_{2}$ is used instead of $\mathbb{R}^N$, and so $N=\infty$.  
\end{Remark}
We now show the proposed kernel-based approximant, i.e.,
\begin{align}\label{Kernel-Based}
    \hat{f}(\bm{x})=\displaystyle \sum_{j=1}^{N}K(\bm{x},\bm{x}_{j})\alpha_{j},\quad \bm{x} \in \mathbb{S}^2,
\end{align}
is an interpolant to a function $f$, where $X=\{\bm{x}_{1},...,\bm{x}_{N}\}$ on $\mathbb{S}^2$ and the values $f(\bm{x}_{1}),...,f(\bm{x}_{N})$ are given.
\begin{proposition}
Assume that for any distinct points $\bm{x}_{1},\ldots \bm{x}_{N}$ on $\mathbb{S}^{2}$, the matrix $A$ has no loss rank.
  The approximation (\ref{Kernel-Based}) is interpolant for given data $\{(\bm{x}_i,f(\bm{x}_i))\}_{i=1}^{N}$.
    \begin{proof}
        According to (\ref{New-Kernel}), for any $\bm{x} \in \mathbb{S}^2$, we can write
    \begin{align}\label{approximate-LS-SVR-Dual-interpolant}
        \hat{f}(\bm{x})
        &=\bm{k}^{T}(\bm{x})\bm{\alpha} \nonumber \\
        &=\bm{k}^{T}(\bm{x})K^{-1}\bm{f} \nonumber \\
        &=\bm{\ell}^{T}(\bm{x})\bm{f},
        \end{align}
        where in the second line of this equality, we have used $K \bm{\alpha}=\bm{f}$, and the fact that $A$ has no loss rank. Due to the definition of $\bm{k}^{T}(\bm{x})$ and
         $\bm{k}^{T}(\bm{x})K^{-1}:=\bm{\ell}^{T}(\bm{x})=[\ell_{1}(\bm{x})\ldots \ell_{N}(\bm{x})]$, it is clear that $\ell_{j}(\bm{x}_{i})=1$ if $i=j$, and $\ell_{j}(\bm{x}_{i})=0$ if $i\neq j$,
which gives the desired results.
    \end{proof}
\end{proposition}

\begin{theorem}
    The proposed kernel based interpolant $\hat{f}$ to a function $f$ in $H$ is a minimum norm interpolant. Indeed
    \begin{align}
        ||\hat{f}||_{H}=\min \{||s||_{H}:\,\,s \in H,\,\, \textsl{with}\,\,s(\bm{x}_{i})=f(\bm{x}_{i}),\,\,i=1,...,N\}.
    \end{align}
    \begin{proof}
        The proof follows from the proof of \cite[Theorem $18.1$]{Fasshauer}.
    \end{proof}
\end{theorem}
\subsection{Smoothing in RKHS}\label{Subsection-Smoothing}
We start with the following definition similar to \cite{Wendland-2006}.
\begin{definition}
Suppose that the RKHS (\ref{Hilbert-Space-1}) is given. For given data $X=\{\bm{x}_{1},...,\bm{x}_{N}\}$ on $\mathbb{S}^2$ and function values $\bm{f}=\{f(\bm{x}_{1}),...,f(\bm{x}_{N})\}$, where $f(\bm{x}_{i}) \in \mathbb{R},\,\,i=1,...,N$, we define 
 \begin{align}\label{Smoothing-Functional}
 F(\hat{f}):=\displaystyle\sum_{j=1}^{N}(f(\bm{x}_{j})-\hat{f}(\bm{x}_{j}))^2+\dfrac{1}{\gamma}||\hat{f}||^2_{H},
 \end{align}
 where $\frac{1}{\gamma}$ is smoothing parameter. Then, we find the solutions $\hat{f}_{\gamma} \in H$ to
 \begin{align}\label{Min-Main}
     \min_{\hat{f} \in H}F(\hat{f}).
 \end{align}
\end{definition}
The following lemma shows that the unique solution of (\ref{Min-Main}) is the same with problem (\ref{LS-SVM-1})-(\ref{LS-SVM-2}).
\begin{lemma}\label{Lemma-Smoothing-Functional}
    The unique solution of problem (\ref{Min-Main}) with target function $f\in H$ is $\hat{f}_{\gamma}$ given by
    \begin{align}
        \hat{f}_{\gamma}=\displaystyle \sum_{j=1}^{N}K(.,\bm{x}_{j})\alpha_{j},
    \end{align}
    where $\bm{\alpha}$ can be computed by solving the following linear system.
    \begin{align}
        (K+\frac{1}{\gamma}I_{N})\bm{\alpha}=\bm{f}.
    \end{align}
    \begin{proof}
First, it should be noted that $K$ is assumed to be  positive definite. So, the proposed functional (\ref{Smoothing-Functional}), which can be written with respect to the vector $\bm{\alpha}$, i.e., $F(\bm{\alpha})=:\dfrac{1}{\gamma}\bm{\alpha}^{T}K\bm{\alpha}+(K\bm{\alpha}-\bm{f})^{T}(K\bm{\alpha}-\bm{f})$, is convex. Hence, it has a unique minimum. By differentiating the functional with respect to $\bm{\alpha}$, we have
\begin{align}\label{LS-SVM-Regular-2}
\nabla_{\bm{\alpha}}F(\bm{\alpha})=\dfrac{2}{\gamma}K\bm{\alpha}+2 K^{T}(K\bm{\alpha}-\bm{f}).
\end{align}
Considering $\nabla_{\bm{\alpha}}F(\bm{\alpha})=\bm{0}$ (the necessary and sufficient condition for finding a unique minimum), we get
\begin{align}\label{LS-SVM-Regular-3}
(\dfrac{2}{\gamma}K+2K^{T}K)\bm{\alpha}=2K^{T}\bm{f}.
\end{align}
Since $K$ is a symmetric positive definite matrix, multiplying both sides of (\ref{LS-SVM-Regular-3}) by the left inverse $K$ and dividing them to factor $2$, we achieve
\begin{align}\label{LS-SVM-Regular-4}
(K+\dfrac{1}{\gamma} I_{N})\bm{\alpha}=\bm{f},
\end{align}
which completes the proof.
\end{proof}
\end{lemma}
We now give the following proposition showing an upper bound for the approximation $\hat{f}_{\gamma}$ at points $X$. 
\begin{proposition}\label{Propposition-Main}
    Suppose $\hat{f}_{\gamma} \in H$  is a solution of \eqref{Min-Main}. 
    Also, assume that $f(\bm{x}_i)$ for $i=1,...,N$ are given, where $f\in H$. Besides, suppose that there exists a function
$\hat{f}_\infty \in H$ with $\displaystyle\sum_{j=1}^{N}(f(\bm{x}_{j})-\hat{f}_{\infty}(\bm{x}_{j}))^2=0$ and $||\hat{f}_\infty||^2_{H} \leq ||f||^2_{H}$. Then, we have
\begin{align}
&||\hat{f}_{\gamma} ||^2_{H} \leq ||f||^2_{H}\nonumber \\
&|f(\bm{x}_{i})-\hat{f}_{\gamma}(\bm{x}_{i})|\leq \gamma^{-1/2}||f||_{H},\quad i=1,...,N.
\end{align}
\begin{proof}
    For the proof, see \cite[Proposition $3.1$]{Wendland-2006}.
\end{proof}
\end{proposition}
In kernel-based methods, it is most common to represent the above upper bounds in terms of the fill distance of $X$. Hence, if we set $\gamma \geq h^{-\beta
}_{X,\mathbb{S}^2}$, where $\beta$
can be chosen to deal with a desired convergence rate, then the above error bounds at points $X$ obey
\begin{align}\label{Upper-Bound-Error}
|f(\bm{x}_{i})-\hat{f}_{\gamma}(\bm{x}_{i})|\leq h^{\beta/2
}_{X,\mathbb{S}^2}||f||_{H},\quad i=1,...,N.
\end{align}
Note that, to deal with a proper rate of convergence, it is recommended to set $\beta \geq 2$ that is an appropriate way to choose parameter $\gamma$. Of course, there are other algorithms to select this parameter such as generalized cross validation (cf. \cite{Fasshauer-Kernel}).


\subsection{Connection to the KRR}

We can utilize the kernel-based approximant  (\ref{approximate-LS-SVR-Dual}) to solve the main problem (\ref{LS-SVR-1}) via the collocation technique directly. Hence, we have
\begin{align}\label{LS-SVR-Dual-Approximation-1}
\sum_{j=1}^{N}\mathcal{L}k(\bm{x}_{i},\bm{x}_{j})\alpha_{j}=f(\bm{x}_{i}),\quad i=1,\ldots N,
\end{align}
in which $\mathcal{L}k(\bm{x}_{i},\bm{x}_{j})=\langle \mathcal{L}\bm{a}(\bm{x}_{i}),\mathcal{L}\bm{a}(\bm{x}_{j})\rangle,\, i,j=1,\ldots N$. So, Eq. (\ref{LS-SVR-Dual-Approximation-1}) can be written as:
\begin{align}\label{LS-SVR-Dual-Approximation-2}
AA^{T}\bm{\alpha}=\bm{b}.
\end{align}
Now, assume that $A$ has full rank.  Let $\bm{w}^{*}$ be the exact solution to the linear system $A\bm{z}=\bm{b}$. So, 
the linear system (\ref{LS-SVR-Dual-Approximation-2}) is  known as the normal equations of the following least squares problem (cf. \cite[Chapter $8$]{Saad}):
\begin{align}\label{LS-SVR-Dual-Normal}
\min ||\bm{w}^{*}-A^{T}\bm{\alpha}||_{2}.
\end{align}
By defining $A^{T}\bm{\alpha}=\bm{z}$,  problem (\ref{LS-SVR-Dual-Normal}) finds the solution vector $\bm{z}$ that is closest to $\bm{w}^{*}$ in two-norm sense. 

\begin{lemma}
	Assume that $A$ has full rank. The solution of problem (\ref{LS-SVR-Dual-Normal}) is given by:
	\begin{align}
		\bm{w}=A^{T}\bm{\alpha}^{*},
	\end{align}
	where
	\begin{align}\label{LS-SVR-Normal-1}
	\bm{\alpha}^{*}=\displaystyle \sum_{i=1}^{N}\dfrac{\bm{v}^{T}_{i}\bm{b}}{\sigma^{2}_{i}}\bm{v}_{i}.
	\end{align}
	in which $\bm{v}_{i}$ is the $i$-th column of the orthogonal matrix $V$ from the singular value decomposition (SVD) of matrix $A^T$. Besides, for $i=1,...,N$, $\sigma_{i}$ are the singular values of $A^T$.
	\begin{proof}
		As mentioned earlier, assume that $\bm{w}^{*}$ is the exact solution to the linear system $A\bm{z}=\bm{b}$.
		Using the SVD of matrix $A^T$, we can write
		\begin{align}\label{LS-SVR-Normal-2}
		\left\| {{\bm{w}^{*}} - {A^T}{\bm{\alpha }}} \right\|_2^2 &= \left\| {{\bm{w}^{*}} - \left( {U\sum {V^T}} \right){\bm{\alpha }}} \right\|_2^2 \nonumber \\
		&= \left\| {{U^T}\left( {{\bm{w}^{*}} - \left( {U\sum {V^T}} \right){\bm{\alpha }}} \right)} \right\|_2^2 \nonumber \\
		&= \left\| {{U^T}{\bm{w}^{*}} - \sum \underbrace {\left( {{V^T}{\bm{\alpha }}} \right)}_{:=\bm{y}}} \right\|_2^2 \nonumber \\
		& = \displaystyle\sum\limits_{i = 1}^N {{{\left( {{\sigma _i}{y_i} - \left( {{\bm{u}}_i^T{\bm{w}^{*}}} \right)} \right)}^2}.} 
		\end{align} 
		Accordingly, to deal with the minimum of problem (\ref{LS-SVR-Dual-Normal}), we should have
		\begin{align}\label{LS-SVR-Normal-3}
		y_{i}=\dfrac{{{\bm{u}}_i^T{\bm{w}^{*}}}}{\sigma_{i}},\quad i=1,...,N.
		\end{align}
		Since $\bm{y}=V^{T}\bm{\alpha}$, we get $\bm{\alpha}^{*}=V\bm{y}$. And, considering relation (\ref{LS-SVR-Normal-2}) gives
		\begin{align}\label{LS-SVR-Normal-3-1}
		\bm{\alpha}^{*}=\displaystyle \sum_{i=1}^{N}\dfrac{\bm{u}^T_{i}\bm{w}^{*}}{\sigma_{i}}\bm{v}_{i}.
		\end{align}
		On the other hand,  
		applying the SVD to $A$, we have $V \Sigma^{T}U^{T} \bm{w}^{*}=\bm{b}$, which yields $\Sigma^{T}U^{T} \bm{w}^{*}=V^{T}\bm{b}$. Finally, we get 
		$U^{T}\bm{w}^{*}=(\Sigma^{T})^{-1}V^{T}\bm{b}$, and replacing it into (\ref{LS-SVR-Normal-3-1}), the desired result is concluded.
	\end{proof} 
\end{lemma}

On the other hand, it is possible to find an upper bound for $\left\| {{\bm{w}^{*}} - \bm{w}} \right\|_2$ in terms of fill distance of points $X$ based on Schaback's theory \cite{Schaback} and some details in \cite{mirzaei:2017-1}. To do so,  
assume that $\bm{w}^{*}$ is the vector of exact values so that $A\bm{w}^{*}=\bm{b}$. If $\bm{w}^{*}$ is replaced by the vector of approximate values, $\tilde{\bm{w}}$, then according to problem (4.24), we have
\begin{equation}\label{Main-error-1}
||\bm{w}^{*}-\bm{w}||_{2}=\min ||\bm{w}^{*}-\bm{z}||_{2} \leq ||\bm{w}^{*}-\tilde{\bm{w}}||_{2}. 
\end{equation}
Now, we can find the upper bound for $||\bm{w}^{*}-\bm{w}||_{2}$ based on  \cite{mirzaei:2017-1}, yielding
\begin{equation}\label{Main-error-2}
||\bm{w}^{*}-\tilde{\bm{w}}||_{2} \leq C C_{s}(A)(1+K(A))h^{m+1}_{X,\mathbb{S}^2}||{w}^{*}||_{C^{m+1}(\mathbb{S}^2)},
\end{equation}
in which $C$ is a positive constant, $C_{s}(A)$ is the stability constant, $K(A)$ is equal to $1$ if $\tilde{\bm{w}}$ is calculated via the minimization of the residual $||A\bm{w}^{*}-\bm{b}||_{2}$, or to $0$ if $\bm{b}$ is in the range of $A$.

Hence, we can finally find the following upper bound.
\begin{equation}\label{Main-error-3}||\bm{w}^{*}-{\bm{w}}||_{2} \leq C C_{s}(A)(1+K(A))h^{m+1}_{X,\mathbb{S}^2}||{w}^{*}||_{C^{m+1}(\mathbb{S}^2)}.
	\end{equation}

\begin{theorem}\label{equivalent-opti-problem}
Under the above assumptions for matrix $K$, the regularized least squares approximation problem, i.e.,
\begin{align}\label{LS-SVM-Regular-1}
\min_{\bm{\alpha}} \dfrac{1}{2}\bm{\alpha}^{T}K\bm{\alpha}+\dfrac{\gamma}{2}(K\bm{\alpha}-\bm{b})^{T}(K\bm{\alpha}-\bm{b}).
\end{align}
has the same solution with the problem (\ref{LS-SVM-1})-(\ref{LS-SVM-2}).
\begin{proof}
The proof can be done in the same manner presented in Lemma \ref{Lemma-Smoothing-Functional}.
\end{proof}
\end{theorem}

\begin{corollary}
	Assume that $\bm{\alpha}$ is the solution of linear system (\ref{LS-SVM-Regular-4}). If $\gamma \to +\infty$, then $\bm{\alpha} \to \bm{\alpha}^{*}$, where $\bm{\alpha}^{*}$ is obtained as (\ref{LS-SVR-Normal-1}).
	\begin{proof}
		If $A^{T}=U\Sigma V^{T}$ is the SVD of $A^{T}$, then (\ref{LS-SVM-Regular-4}) changes to
		\begin{align}
\begin{array}{l}
\left( {V{\sum ^T}{U^T}U\sum {V^T} + \frac{1}{\gamma }{I_N}} \right){\bm{\alpha }} = {\bm{b}}\overbrace  \Rightarrow ^{{U^T}U = {I_N}}\left( {V{\sum ^T}\sum {V^T} + \frac{1}{\gamma }{I_N}} \right){\bm{\alpha }} = {\bm{b}}\\ \\
\Rightarrow \left( {{V^T}V{\sum ^T}\sum {V^T} + \frac{1}{\gamma }{V^T}{I_N}} \right){\bm{\alpha }} = {V^T}{\bm{b}}\overbrace  \Rightarrow ^{{V^T}V = {I_N}}\left( {{\sum ^T}\sum  + \frac{1}{\gamma }{I_N}} \right)\left( {{V^T}{\bm{\alpha }}} \right) = {V^T}{\bm{b}}.
\end{array}
		\end{align}
		Hence, we have
		\begin{align}
		\bm{\alpha}=\displaystyle \sum_{i=1}^{N}\dfrac{\bm{v}^{T}_{i}\bm{b}}{\sigma^{2}_{i}+\frac{1}{\gamma}}\bm{v}_{i}.
		\end{align}
		Now, if $\gamma \to +\infty$, the desired result is obtained.
	\end{proof}
\end{corollary}

\begin{Remark}
    The problem (\ref{LS-SVM-Regular-1}) is known as ridge regression, and it can be observed as the diagonal stabilization/regularization of the usual approximation linear system (\ref{LS-SVR-Dual-Approximation-2}). 
\end{Remark}

\begin{proposition}
    Assume that $\mathcal{L}=\texttt{Id}$ in problem (\ref{LS-SVR-1}). Suppose also that all assumptions proposed in Theorem \ref{equivalent-opti-problem} are satisfied. If $\gamma =+\infty$, then the approximate solution (\ref{approximate-dual-LS-SVR}) is interpolant for given data $\{(\bm{x}_i,f(\bm{x}_i))\}_{i=1}^{N}$ (see, e.g., \cite{Fasshauer, Fasshauer-Kernel}).
    \begin{proof}
        According to (\ref{approximate-LS-SVR-Dual}) and (\ref{LS-SVM-Regular-4}), for any $\bm{x} \in \mathbb{S}^2$, we can write
    \begin{align}\label{approximate-LS-SVR-Dual-interpolant-1}
        \hat{f}_{\gamma}(\bm{x})
        &=\bm{k}^{T}(\bm{x})(K+\dfrac{1}{\gamma} I_{N})^{-1}\bm{b} \nonumber \\
        &=\bm{k}^{T}(\bm{x})K^{-1}(K+\dfrac{1}{\gamma} I_{N})^{-1}K\bm{b}\nonumber \\
        &=\bm{\ell}^{T}(\bm{x})(K+\dfrac{1}{\gamma} I_{N})^{-1}K\bm{b},
        \end{align}
        where in the second line of this equality, we have used the push-through identity (cf. \cite{Fasshauer-Kernel}), and $\bm{k}^{T}(\bm{x})K^{-1}:=\bm{\ell}^{T}(\bm{x})=[\ell_{1}(\bm{x})\ldots \ell_{N}(\bm{x})]$. Due to the definiton of $\bm{k}^{T}(\bm{x})$, it is clear that $\ell_{j}(\bm{x}_{i})=1$ if $i=j$, and $\ell_{j}(\bm{x}_{i})=0$ if $i\neq j$. Considering $\gamma=+\infty$ and using the fact that $\bm{b}=[f(\bm{x}_{i})]_{1 \leq i \leq N}$, we thus have
        \begin{align}
        \nonumber \\
        \hat{f}_{\infty}(\bm{x})=\bm{\ell}^{T}(\bm{x})\bm{b},
        \end{align}
which gives the desired results.
    \end{proof}
\end{proposition}

\begin{Remark}
	If we assume that $f(\bm{x}_{i})$ is contaminated by error $e_{i}$, i.e., $f(\bm{x}_{i})+e_{i}$ for $i=1,...,N$, then the new righ hand side is given by $\bm{b}=[f(\bm{x}_{i})+e_{i}]_{1\leq i \leq N}$. These can be viewed as the noisy data. According to the approximation (\ref{approximate-LS-SVR-Dual-interpolant}), the term 
	$\tilde{\bm{b}}=:(K+\dfrac{1}{\gamma} I_{N})^{-1}K\bm{b}$ refers to the smoothed data or filtered data. Hence, relation (\ref{approximate-LS-SVR-Dual-interpolant}) can be known as the kernel interpolation for these filtered data. In this case, we have
	    \begin{align}\label{approximate-LS-SVR-Dual-interpolant-2}
	\hat{f}(\bm{x})=\bm{\ell}^{T}(\bm{x})\tilde{\bm{b}},
	\end{align}
\end{Remark}

\section{Application to two phase field models on $\mathbb{S}^2$}\label{sect-LS-SVR-AC-CAC}
In this section, we derive the new fully discrete schemes for both
phase-field models, i.e., Eqs. (\ref{e1}) and (\ref{e1-CAC}). We begin by applying an implicit
time-splitting approach of accuracy order one to discretize each equation in time.  Then, we develop the presented LS-SVR for the spatial discretization of the obtained semi-discrete equations.

\subsection{An implicit time-splitting
technique}
The sequential time-splitting idea
\cite{xiao-et-al:2017} on Eq. \eqref{e1} and a type of first-order operator splitting approach \cite{Kim11} on \eqref{e1-CAC} have been applied (see also
\cite{MohDehMir}).

To start, suppose that the time interval is $[0, T]$, where $T$ is
the final time. To partition this interval, consider a uniform time step
$\Delta t=\frac{T}{K}$ with the number of sub-intervals to be $K$.
So, $[0,T]=\bigcup_{n=1}^{K}[t_{n-1},t_{n}]$, $t_{0}=0$, $t_{n}=n
\Delta t$ for $n=1,...,K-1$ and $t_{K}=K\Delta t=T$. Hence, the proposed sequential-implicit splitting method could be executed for
Eq. \eqref{e1} as follows:
\begin{align}
(\texttt{id}-\Delta t \texttt{L}) u^{n^{*}}&=u^{n}, \label{e5} \\
u^{n+1}&=\frac{u^{n^{*}}}{\sqrt{\exp\left({-\frac{2 \Delta
t}{{\epsilon}^{2}}}\right)+ (u^{n^{*}})^{2}\left(1-\exp\left({-\frac{2
\Delta t}{{\epsilon}^{2}}}\right)\right)}},\quad n=0,\ldots,K-1.
\label{e6}
\end{align}
 In the above
semi-discrete equations, $u^n$ shows the approximate value of
the exact solution $u(\cdot,t_n)$,
 and $u^{n^{*}}$ is the approximate intermediate solution.

Also, by applying a type of time-splitting technique utilized in \cite{Kim11}
to Eq. (\ref{e1-CAC}), we get the following semi-discrete equations:
\begin{align}
(\texttt{id}-\Delta t \texttt{L}) u^{n^{*}}&=u^{n}, \label{e5-CAC} \\
u^{n+1^{**}}&=\frac{u^{n^{*}}}{\sqrt{\exp\left({-\frac{2 \Delta
t}{{\epsilon}^{2}}}\right)+ (u^{n^{*}})^{2}\left(1-\exp\left({-\frac{2
\Delta t}{{\epsilon}^{2}}}\right)\right)}},
\label{e6-CAC} \\
u^{n+1}&=u^{n+1^{**}}+\Delta t\beta^{n+1^{**}} \sqrt{F(u^{n+1^{**}})}, \label{e7-CAC}
\quad n=0,\ldots,K-1,
\end{align}
in which $u^{n+1^{**}}$ is known as the second approximate intermediate solution, and
\begin{align}\label{e8-CAC}
    \beta^{n+1^{**}}=\dfrac{\displaystyle \int_{\mathbb{S}^2}(u^{0}-u^{n+1^{**}})
    ~d\bm{x}}{ \Delta t \displaystyle \int_{\mathbb{S}^2}\sqrt{F(u^{n+1^{**}})}~d\bm{x}}.
\end{align}
In Appendix A, we have proved that the proposed time discretization (\ref{e5-CAC})-(\ref{e7-CAC}) preserves the total mass in time.

\subsection{The LS-SVR for AC and CAC models}
Here, we show how the presented LS-SVR can be applied to the
semi-discrete problems (\ref{e5})-(\ref{e6}) and
(\ref{e5-CAC})-(\ref{e7-CAC}). For doing so,
the set of $N $ scattered trial points such as
$X=\{\bm{x}_{1},...,\bm{x}_{N}\}$ on $\mathbb{S}^2$ is considered. Hence, the approximate
solution at each time step $n$ can be written as:
\begin{align}\label{LS-SVM-4}
    \hat{u}(\bm{x},t^{n})=\bm{a}^{T}(\bm{x})\bm{w}_{n},\quad \bm{x} \in \mathbb{S}^2,
\end{align}
where $\bm{w}_{n}=[w^n_{1},...,w^n_{N}]^{T}$, and
$\bm{a}^{T}(\bm{x})=[a_{1}(\bm{x}),...,a_{N}(\bm{x})]$ denotes the vector with trial functions.

\subsubsection{AC model}
We first note that the
linear differential operator in both studied mathematical models is
$\mathcal{L}=\texttt{Id}-\Delta t \texttt{L}$.
Replacing the approximate solution
(\ref{LS-SVM-4}) and executing the collocation technique based on the GMLS approximation, we get the following problem per time step
$n=0,..., K-1$;

Solve the following minimization problem.
\begin{align}
    &\min_{\bm{w}^{*}_n,\bm{e}^{*}_n \in  \mathbb{R}^{N}} \dfrac{1}{2}\bm{w}^{*^{T}}_n\bm{w}^{*}_n+\dfrac{\gamma}{2}\bm{e}^{*^T}_n\bm{e}^{*}_n,\label{LS-SVM-5}\\
    &s.t.~A_{X}\bm{w}^{*}_{n}-\Delta t L_X\bm{w}^{*}_{n}=A_{X}U^{n}_{X}+\bm{e}^{*}_{n},\label{LS-SVM-6}
\end{align}
where $\bm{e}_{n}$ is the residual term,
$A_{X}=[a_{j}(\bm{x}_{i})]_{1\leq i,j \leq N}$ and $L_{X}=[a^{\texttt{L}}_{j}(\bm{x}_{i})]_{1\leq i,j \leq N}$ are sparse matrices thanks to the presented GMLS approximation, which are the same with those differentiation matrices computed for approximating the point-evaluation functional and Laplace-Beltrami operators in Section \ref{GMLS-spheres}, respectively.
 and $\bm{w}^{*}_{n}$ is the vector of intermediate weights, which can be computed by applying KKT conditions as follows:
\begin{align}\label{LS-SVM-7}
    (A_{L}A^{T}_{L}+\frac{1}{\gamma}I_{N})\bm{\alpha}^{*}_{n}=A_{X}U^{n}_{X},
\end{align}
which is a dual form of the problem
(\ref{LS-SVM-5})-(\ref{LS-SVM-6}). Besides, $U^{n}_{X}$ is the
approximate solution at points $X$ in $t=t_{n}$ with
$A_{L}=A_X-\Delta t L_{X}$, and
\begin{align}\label{LS-SVM-7-1}
    \bm{w}^{*}_{n}=A^{T}_{L}\bm{\alpha}^{*}_{n},\quad \bm{e}^{*}_{n}=-\dfrac{1}{\gamma}\bm{\alpha}^{*}_{n}.
\end{align}
We can now compute the approximate solution to the intermediate variable as
\begin{align}\label{LS-SVM-8}
    U^{n*}_X=A_{X}\bm{w}^{*}_{n}.
\end{align}
Finally, due to Eq. (\ref{e6}), the approximate solution at time level $t=t_{n+1}$ can be calculated with the following relation.
\begin{align}\label{LS-SVM-9}
U^{n+1}_X&=\frac{U^{n*}_X}{\sqrt{\exp\left({-\frac{2 \Delta
t}{{\epsilon}^{2}}}\right)+ (U^{n*}_X)^{2}\left(1-\exp\left({-\frac{2
\Delta t}{{\epsilon}^{2}}}\right)\right)}}.
\end{align}
Now, we can give the following theorem that shows the proposed numerical scheme applied to the AC equation (\ref{e1}) has a unique solution.
\begin{theorem}\label{Th-LS-SVR-2}
For $n \geq 0$, the discrete scheme given in (\ref{LS-SVM-7})-(\ref{LS-SVM-9}) has a unique solution $U^{n+1}_{X}$.
\begin{proof}
    For $n \geq 0$, the linear system (\ref{LS-SVM-7}) has a unique solution $\bm{\alpha}^{*}_{n}$ because the coefficient matrix $B_{L}=:A_{L}A^T_{L}+\frac{1}{\gamma}I_{N}$ is positive definite (see the proof of Theorem \ref{Th-LS-SVR-1}). Accordingly, it is invertible, and hence the vectors (\ref{LS-SVM-7-1}) and the approximate solution of the intermediate variable $u^{n^{*}}$, i.e., $U^{n*}_X$ are uniquely obtained. Finally, since $U^{n+1}_{X}$ is directly computed by relation (\ref{LS-SVM-9}), the proof is completed.
\end{proof}
\end{theorem}

\subsubsection{CAC model}
Additionally, to apply the proposed LS-SVR to the semi-discrete problem (\ref{e5-CAC})-(\ref{e6-CAC}), we employ the same procedure executed to the AC model. Accordingly, using
the collocation technique based on the GMLS approximation with $N$
scattered trial points on $\mathbb{S}^{2}$, we have the following.

\begin{figure}[h]
    \centering
    \includegraphics[width=0.49\textwidth]{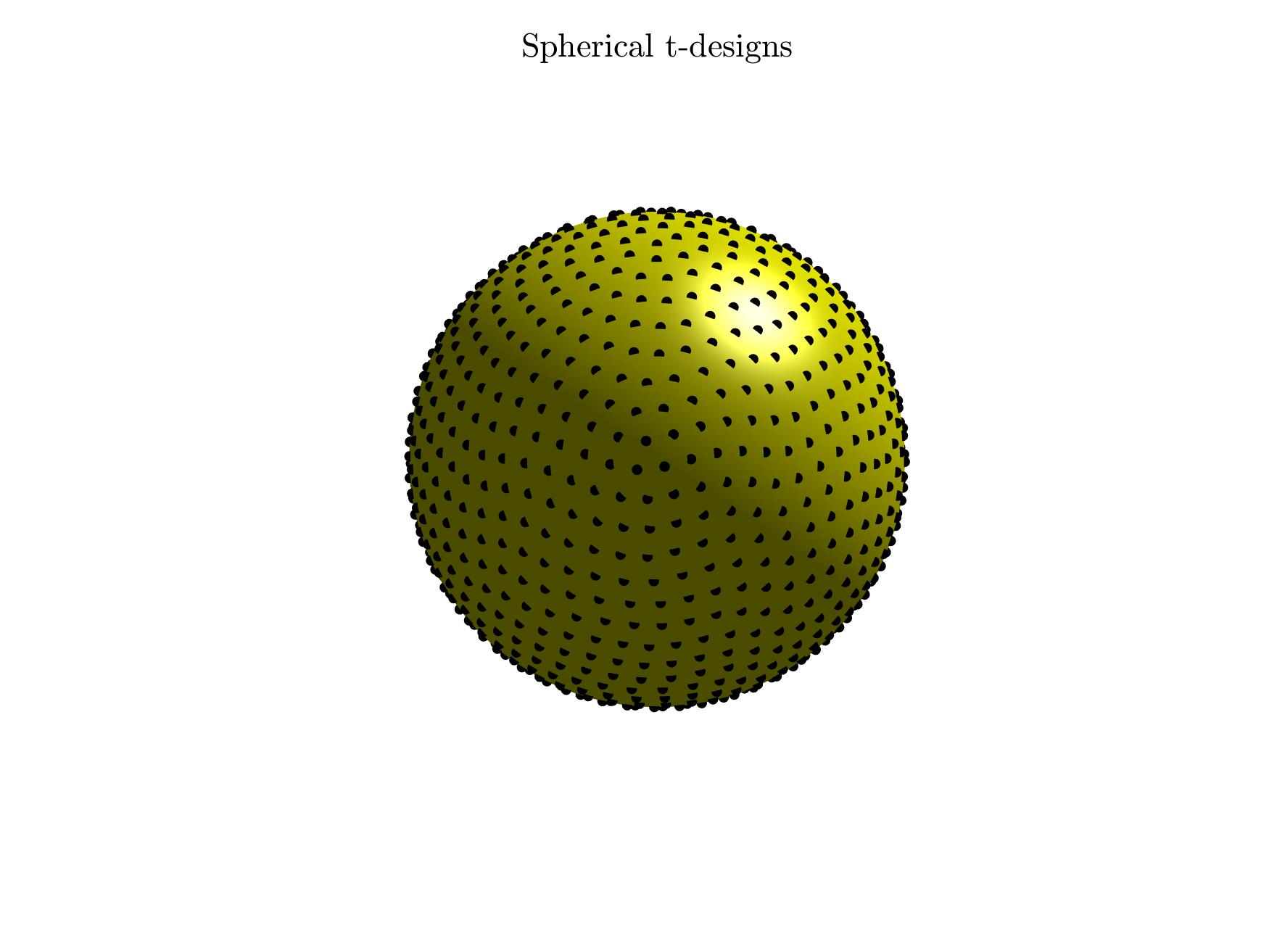}
    \includegraphics[width=0.49\textwidth]{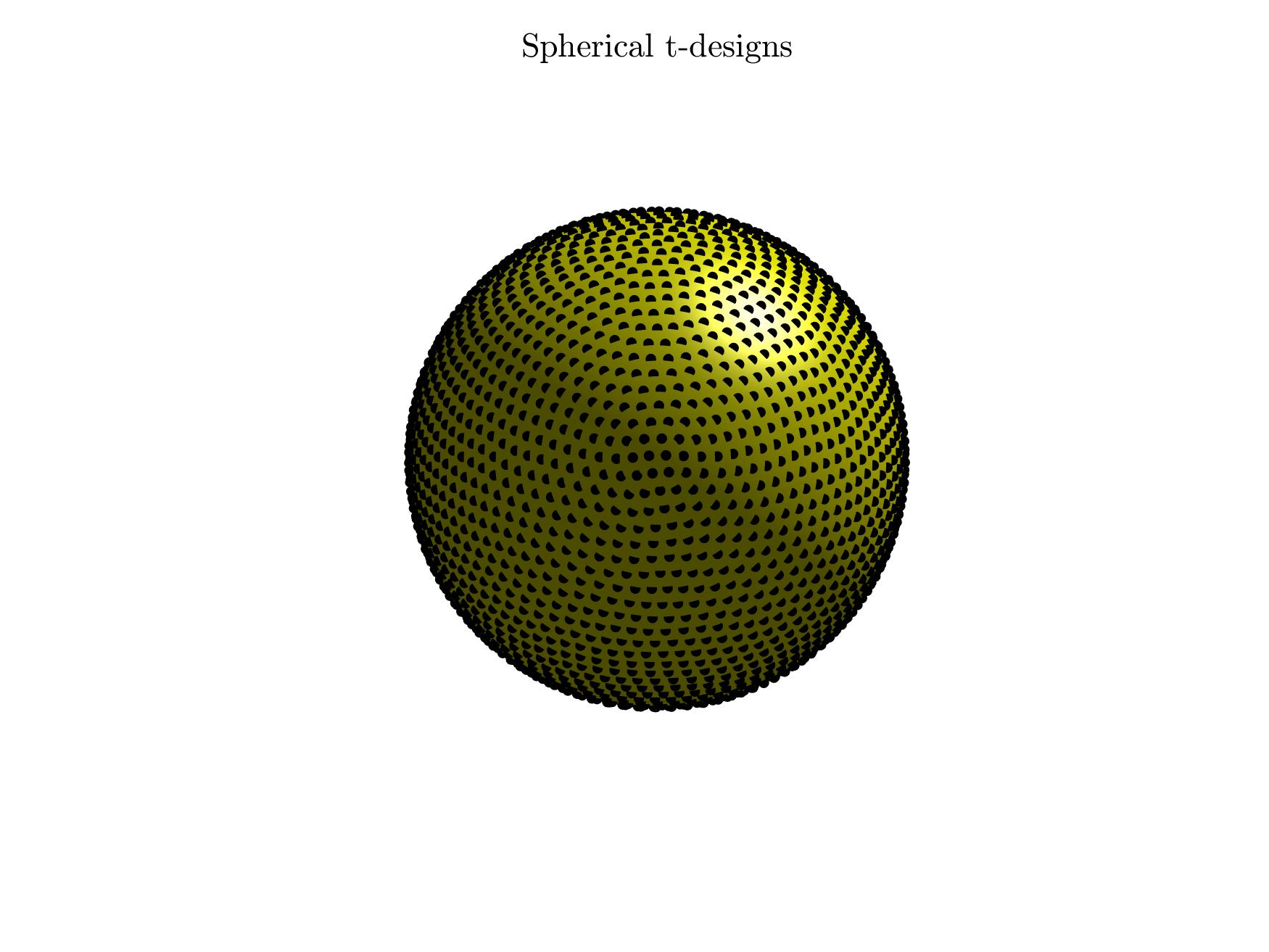}
    \caption{The plots of $966$ (left) and $2666$ (right) spherical
        $t$-design points on $\mathbb{S}^2$.} \label{Sphericalt-design}
\end{figure}

For $n=0,..., K-1$, the following minimization problem should be
solved.
\begin{align}
    &\min_{\bm{w}^{*}_n,\bm{e}^{*}_n \in  \mathbb{R}^{N}} \dfrac{1}{2}\bm{w}^{*^{T}}_n\bm{w}^{*}_n+\dfrac{\gamma}{2}\bm{e}^{*^T}_n\bm{e}^{*}_n,\label{LS-SVM-5-CAC}\\
    &s.t.~A_{X}\bm{w}^{*}_{n}-\Delta t L_X\bm{w}^{*}_{n}=A_{X}U^{n}_{X}+\bm{e}^{*}_{n},\label{LS-SVM-6-CAC}
\end{align}
In addition, the vector of intermediate weights, i.e.,
 $\bm{w}^{*}_{n}$ is computed by using the KKT conditions as follows:
\begin{align}\label{LS-SVM-7-CAC}
    (A_{L}A^{T}_{L}+\frac{1}{\gamma}I_{N})\bm{\alpha}^{*}_{n}=A_{X}U^{n}_{X},
\end{align}
that is a dual form of the problem
(\ref{LS-SVM-5-CAC})-(\ref{LS-SVM-6-CAC}). Also, $U^{n}_{X}$ shows
the approximate solution at points $X$ in $t=t_{n}$ and
$A_{L}=A_X-\Delta t L_{X}$, and
\begin{align}\label{LS-SVM-7-1-CAC}
    \bm{w}^{*}_{n}=A^{T}_{L}\bm{\alpha}^{*}_{n},\quad \bm{e}^{*}_{n}=-\dfrac{1}{\gamma}\bm{\alpha}^{*}_{n}.
\end{align}
 The approximate solution to the intermediate variable can now be computed as
\begin{align}\label{LS-SVM-8-CAC}
    U^{n*}_X=A_{X}\bm{w}^{*}_{n}.
\end{align}
Besides, the approximate solution to the second intermediate variable can be calculated with the
following relation due to Eq. (\ref{e6-CAC}).
\begin{align}\label{LS-SVM-9-CAC}
U^{n+1^{**}}_X&=\frac{U^{n^{*}}_X}{\sqrt{\exp\left({-\frac{2 \Delta
t}{{\epsilon}^{2}}}\right)+ (U^{n^{*}}_X)^{2}\left(1-\exp\left({-\frac{2
\Delta t}{{\epsilon}^{2}}}\right)\right)}}.
\end{align}
Finally, according to (\ref{e8-CAC}), the approximate solution at time level $t=t_{n+1}$ is obtained by
\begin{align}\label{LS-SVM-10-CAC}
U^{n+1}_{X}=U^{n+1^{**}}_X+\Delta t\beta^{n+1^{**}}\sqrt{F(U^{n+1^{**}}_{X})},
\end{align}
where
$\beta^{n+1^{**}}$ is computed based on the integration of scattered data
technique applied to (\ref{e8-CAC}), where $16382$ spherical $t$-design points is used \cite{Atkinson-2,Womersley}. For instance, to compute the integration of $u$ over $\mathbb{S}^2$, we have utilized the following approximation.
\begin{equation}\label{Integration}
    \displaystyle \int_{\mathbb{S}^2}u(\bm{x})d\bm{x} \approx
    \displaystyle \dfrac{4 \pi}{N_{g}}\sum_{j=1}^{N_{g}} u(\bm{y}_{j}),
\end{equation}
in which $\{\bm{y}_{1},...,\bm{y}_{N_{g}}\}$ denotes $N_{g}$ spherical $t$-design points (cf. \cite{Atkinson-2,Womersley}). In Figure
\ref{Sphericalt-design}, two sets of spherical $t$-design points
with cardinalities $969$ and $2666$ are plotted.

\begin{theorem}\label{Th-LS-SVR-3}
    For $n \geq 0$, the discrete scheme given in (\ref{LS-SVM-7-CAC})-(\ref{LS-SVM-10-CAC}) has a unique solution $U^{n+1}_{X}$.
    \begin{proof}
        The proof can be done in the same manner presented in Theorem \ref{Th-LS-SVR-2}.
    \end{proof}
    \end{theorem}

In Appendix B, for the reader's convenience, we present the computational algorithms
associated with the numerical methods used to solve the AC and CAC
equations.

\section{Numerical results}\label{sect-numerical}
To demonstrate the potency of the developed numerical method based on a type of SVM, i.e., the LS-SVR combined with a meshfree approximation on spheres, some numerical simulations are presented in this section.

To use the GMLS approximation, we consider the spherical harmonic functions with
degree $m=2$ or less and the weight function as follows (cf. \cite[Chapter $9$]{Wendland}):
 \begin{equation}\label{Weight-Function}
\phi (r) = \left\{ \begin{array}{l}
{(1 - r)^{10}}(429r^4+450r^3+210r^2+50r + 5),\,\,\,\,\,\,\,\,\,\,\,\,\,\,\,\,\,\,\,\,\,\,0 \le r \le 1,\\
0,\,\,\,\,\,\,\,\,\,\,\,\,\,\,\,\,\,\,\,\,\,\,\,\,\,\,\,\,\,\,\,\,\,\,\,\,\,\,\,\,\,
\,\,\,\,\,\,\,\,\,\,\,\,\,\,\,\,\,\,\,\,\,\,\,\,\,\,\,\,\,\,\,\,\,\,\,\,\,\,\,\,\,\,\,\,\,\,\,\,\,\,\,\,\,\,\,\,\,\,\,\,\,\,\,\,\,\,\,\,\,\,\,\,\,\,\,\,\,\,\,\,\,\,\,\,\,\,\,\,\,\,\,\,\,\,\,\,\,\,\,\,\,r
> 1,
\end{array} \right.
\end{equation}
in which $r={dist(\bm{x},\bm{y})}/{\delta}$ and $\bm{x},\bm{y} \in
\mathbb{S}^2$.

\begin{figure}[h!]
\centering
\includegraphics[width=0.45\textwidth]{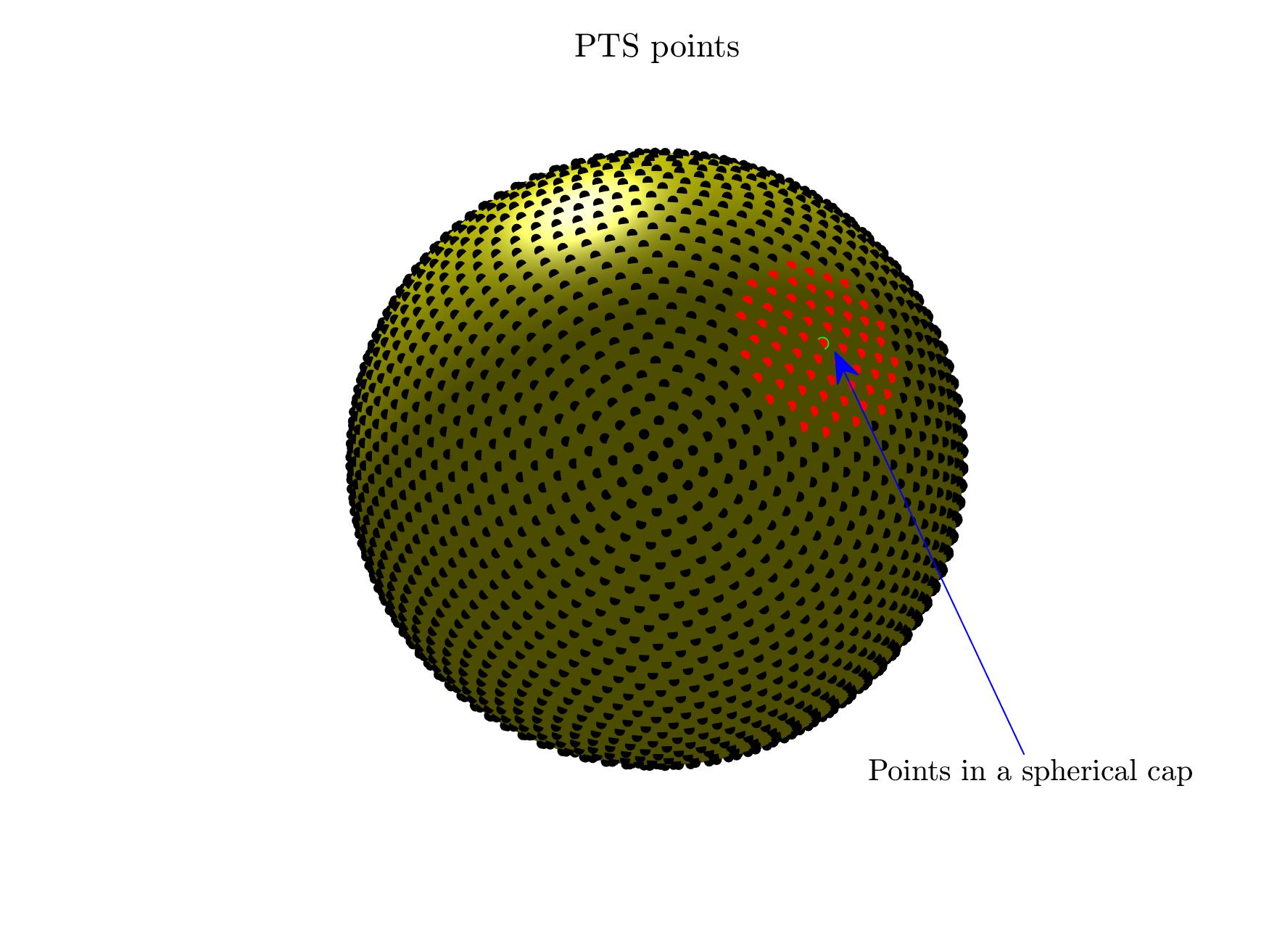}
\includegraphics[width=0.45\textwidth]{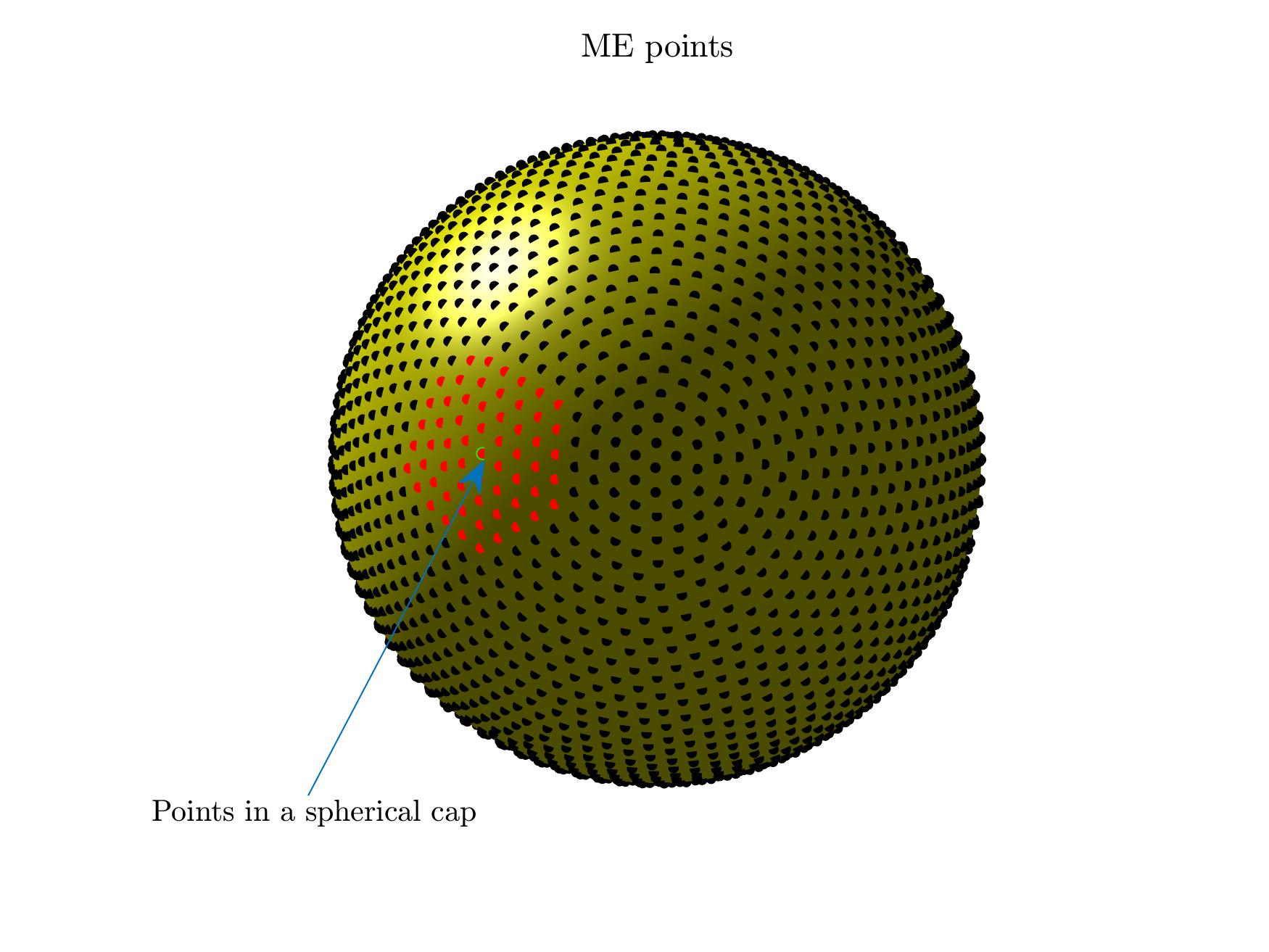}
\caption{The plots of $2601$ ME (left) and PTS (right) points with a
spherical cap containing red points neighboring $2000$-th point.}
\label{points}
\end{figure}

Additionally, the trial approximation (\ref{approximate-LS-SVR}) is constructed by choosing
two sets of quasi-uniform nodes, i.e., the minimum energy (ME) and the
phyllotaxis spiral (PTS) on $\mathbb{S}^2$, where their fill
distance is of order $\mathcal{O}(N^{-1/2})$ (cf.
\cite{Fuselier,MohDehMir}). In Figure \ref{points}, $2601$ ME and
PTS nodes are plotted, where a spherical cap related to a test point and its trial points are denoted.

In the computations presented below, to avoid any instability in solving
the local normal linear systems arising in  Eqs.
(\ref{GMLS-Sphere6}) and (\ref{GMLS-Sphere-new-1}), the QR
factorization of the matrix $\sqrt{W}P$ to compute
$C=WP\big(P^{T}W(\bm{x})P\big)^{-1}$ via the backward substitution
from $RC^{T}=Q^{T}\sqrt{W}$ is applied. Note that, the computational complexity of the approximation \eqref{GMLS-Sphere6} ( or (\ref{GMLS-Sphere-new-1})) for each test point $\bm{x} \in \mathbb{S}^2$ is of order $\mathcal{O}(N(3,m)^3)$ (cf. \cite{mirzaei:2017-1})).

\begin{figure}[h!]
    \centering
    \includegraphics[width=12cm,height=7cm]{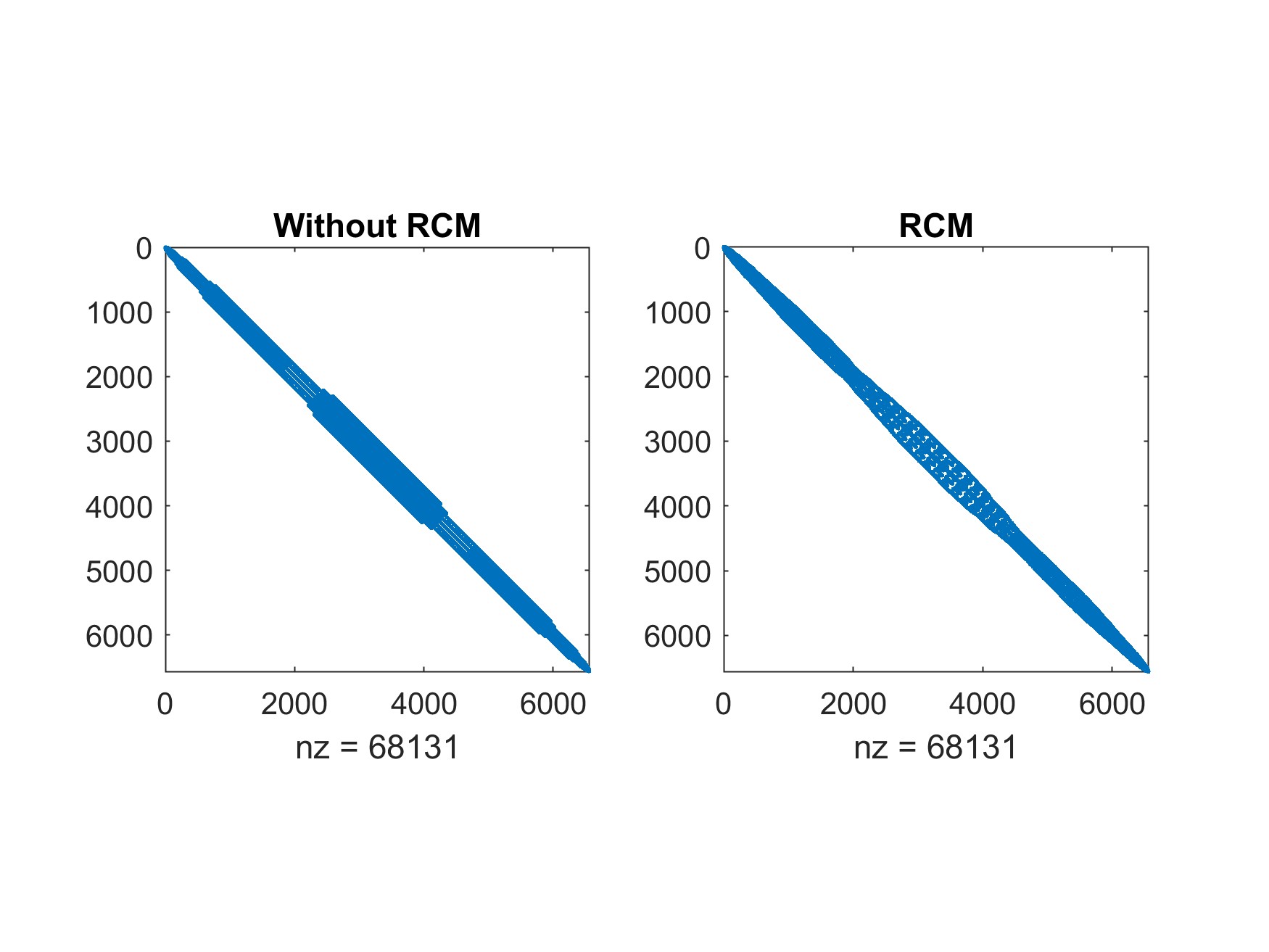}
    \includegraphics[width=12cm,height=7cm]{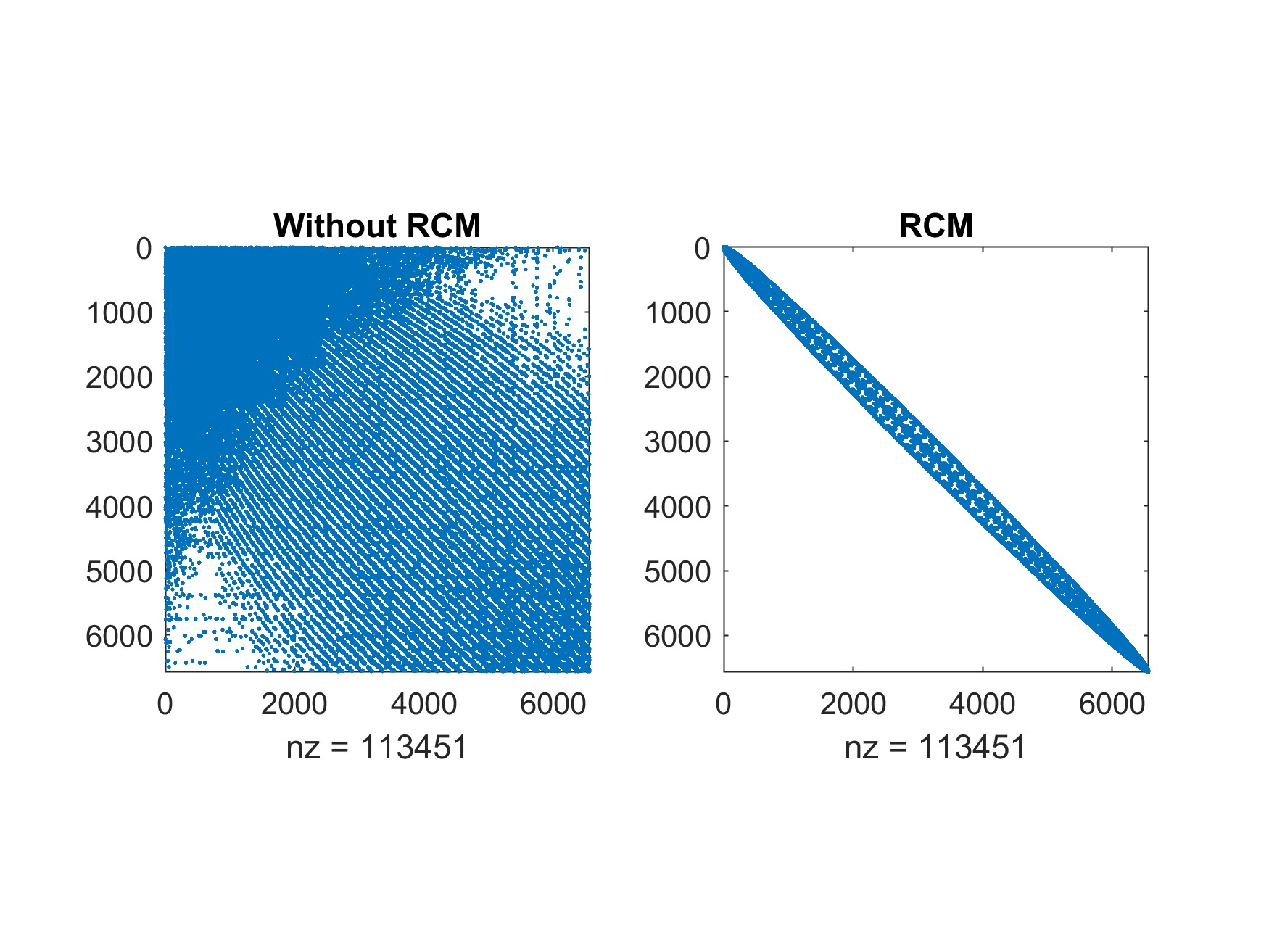}
    \caption{The sparsity pattern of $B_{L}$ before and after re-ordering with the RCM algorithm, where $6561$ PTS points (left) and ME points (right) are used.}
    \label{RCM-Pattern}
\end{figure}

On the other hand, instead of solving the obtained linear systems (\ref{LS-SVM-7}) and (\ref{LS-SVM-7-CAC}), we first use the reverse Cuthill–McKee (RCM) reordering algorithm \cite{George1981,Shankar}.
Then, the preconditioned conjugate gradient (PCG)
algorithm with a zero-fill incomplete lower-upper (ILU(0))
  factorization of the coefficient matrix, obtained by executing the RCM algorithm on $B_{L}$, is employed, where, for stopping the algorithm, tolerance and maximum number of iterations are
chosen $10^{-6}$ and $100$, respectively (cf. \cite{Saad}).  The observed numerical results (not reported here for brevity) show that the ILU(0)  factorization is a proper preconditioner.
Note that the implementation of the RCM approach leads to faster convergence when an iterative solver is used, as it permutes a sparse matrix with a symmetric sparsity pattern into a band matrix form with a small bandwidth (cf. \cite{Shankar}). In Figure \ref{RCM-Pattern}, the sparsity pattern of $B_{L}$ of size $6561 \times 6561$ (assumed $\mathcal{L}=\texttt{Id}$) before and after re-ordering with the RCM algorithm is drawn using both PTS and ME points.


\subsection{Convergence test}\label{Subsec-5.1}
\begin{figure}[h!]
    \centering
    \includegraphics[width=7cm,height=7cm]{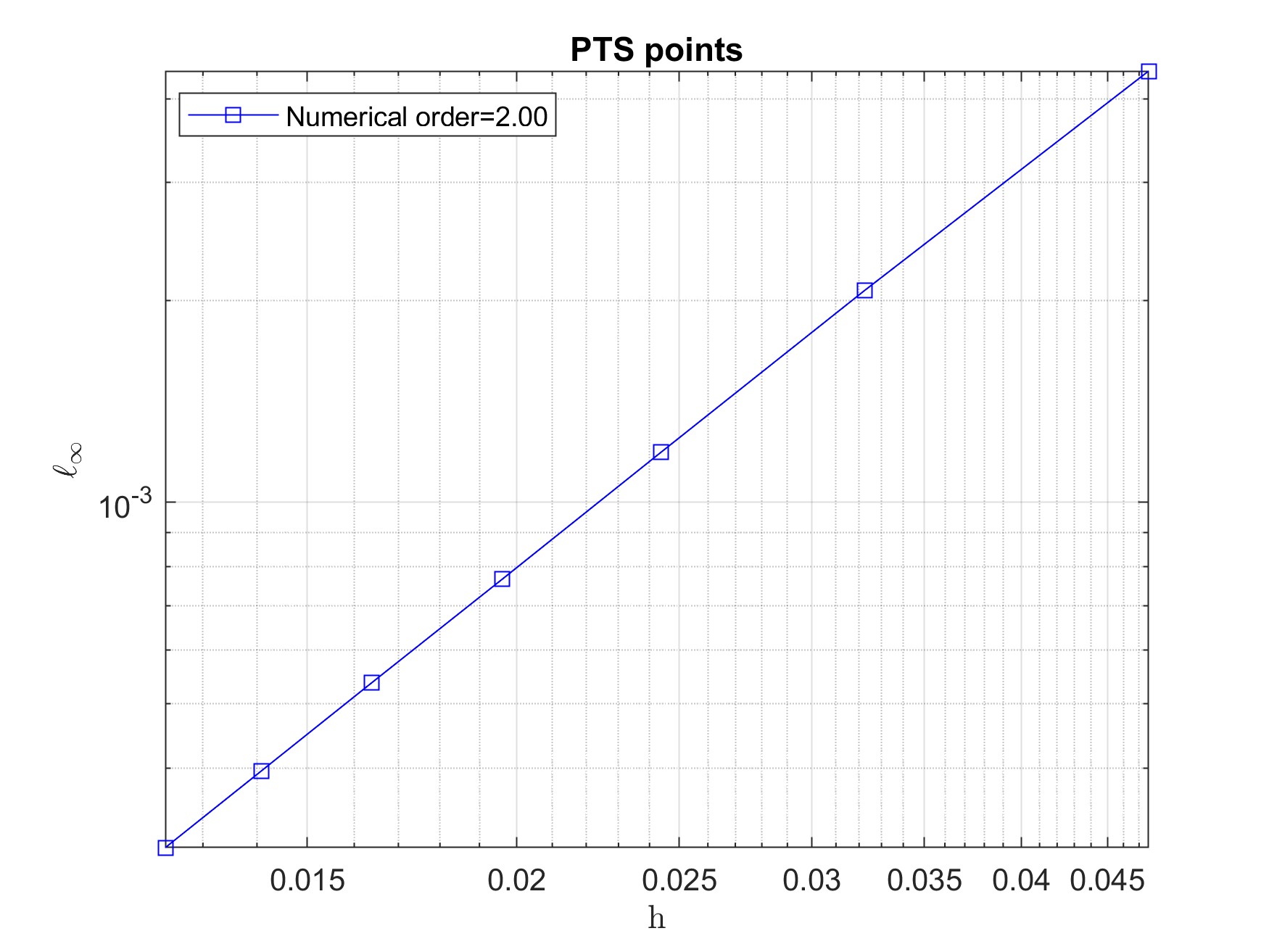}
    \includegraphics[width=7cm,height=7cm]{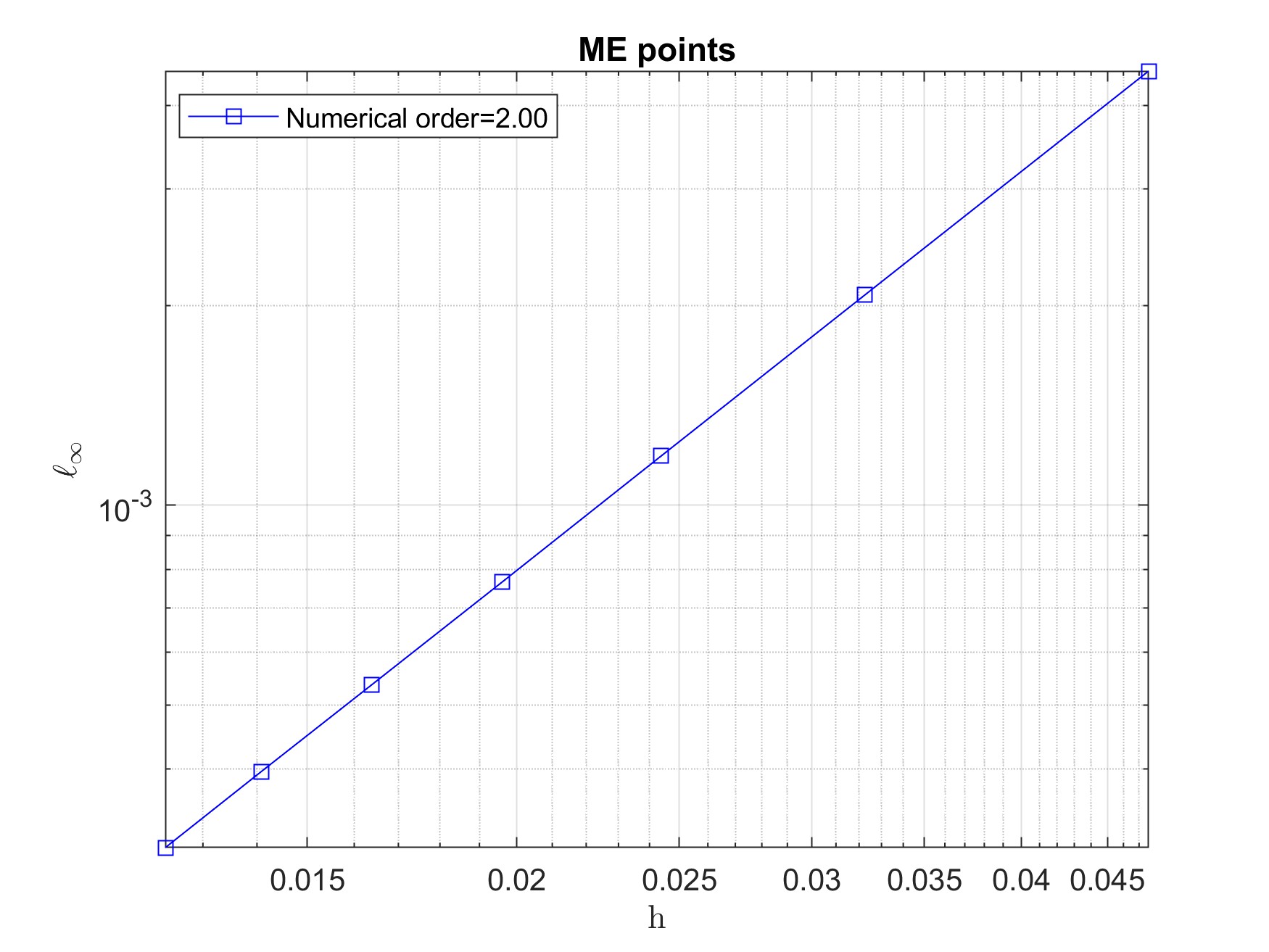}
    \caption{$\ell_{\infty}$ errors with the numerical order, PTS points (left) and ME points (right).}
    \label{fig0-accuracy}
\end{figure}

Here,  we  numerically investigate the error bounds (\ref{Upper-Bound-Error}). To do so, the following smooth test function is considered \cite{Fuselier3,Mir-D-RBF-PU}:
\begin{equation}\label{Exact-Solution-test}
    u(x,y,z)=\cos(10z+3)e^{-x^2}-\sin(4x+5y),\quad (x,y,z) \in \mathbb{S}^2,
\end{equation}
in which $\mathbb{S}^2$ is the unit sphere, i.e.,
\begin{equation}\label{Unit-Sphere}
\mathbb{S}^2=\left\{(x,y,z) \in \mathbb{R}^{3}:x^2+y^2+z^2=1\right\}.
\end{equation}
In all numerical results presented here,
PTS points on $\mathbb{S}^2$ are used.
Besides, the size of each spherical cap is considered to be $\delta=4mh$, where $m=2$ and $h=1/\sqrt{N}$. Besides, $\gamma =h^{-2}$ to deal with the errors, which obey $\mathcal{O}(h)$.
The accuracy of LS-SVR based on the GMLS approximation is measured by computing the discrete maximum errors, i.e., $\ell_{\infty}$ with the numerical order, as displayed in Figure \ref{fig0-accuracy}. However, 
the observed numerical convergence
 demonstrates the results being approximately twice order more accurate than the theoretical rate. From our observations, we also have seen that the numerical orders for different values of $\gamma$, i.e., $h^{-3},h^{-4}$ and $h^{-5}$ are $3$, $4$ and $5$ respectively.

\subsection{Test 1 (Accuracy test): Motion by mean curvature on the unit sphere}

As the first example, we denote the accuracy of numerical scheme in both time and space domains. So,
the AC equation (\ref{e1})
 with the following initial condition is considered \cite{Choi}.
 \begin{equation}\label{initial_motion1}
    u(x,y,z,0)=\tanh\Big(\frac{\sqrt{x^2+y^2}-z}{\sqrt{2}{\epsilon}}\Big),\quad (x,y,z) \in \mathbb{S}^{2},
 \end{equation}
 in which $\epsilon=0.1$.

This problem is known as the motion by mean curvature on the unit sphere (cf. \cite{Choi,MohDehMir}). Indeed, the interface of a circle on the unit sphere is evolving in time due to the motion of surface curvature. Accordingly, let the radii of the sphere and the circle be $R$ and $r$, respectively. By considering some assumptions and simple calculations proposed in \cite{Choi}, the governing equation related to the motion of a circular interface on the sphere is derived as:
\begin{align}\label{Motion}
    \dfrac{dr}{dt}=\dfrac{r^2-R^2}{rR^2},
    \end{align}
where its exact solution is \cite{Choi}
\begin{align}\label{Motion-Exact}
    r(t)=\sqrt{R^2-(R^2-r^2_{0})\exp(\frac{2t}{R^2}))},
\end{align}
where $R=1$ and $r_{0}=\frac{1}{\sqrt{2}}$ \cite{Choi}.

To measure the accuracy of the developed numerical approach in both time and spatial discretizations, the radius of the spherical cap for each test point is chosen to be $\delta=6m h$, where $m=2$ and $h = 1/\sqrt{N}$. Also, the coefficient constant $\gamma$ is empirically chosen to be $10^{12}$.

In Tables \ref{Table-1} and \ref{Table-2}, the discrete maximum errors, i.e., $\ell_{\infty}$ corresponding to the time and space discretizations for approximating the radii (\ref{Motion-Exact}) are reported at $T=0.2$. Note that, in obtaining the errors related to the time discretization, we have used $10201$ PTS points. And, the errors computed based on the spatial discretization are reported by choosing $\Delta t=10^{-4}$. Note also that the average of last three errors is approximately $1.35$ in time, and its value is about $2.93$ in space.

\begin{table}[h!]
    \centering
    \begin{tabular}{ccccccccccccccccccc}
        \hline
        \hline
        $\Delta t$ && $\ell_{\infty}$ && Orders \\
        \hline
        $10^{-2}$ &&$6.99 \times 10^{-2}$&& $-$ \\
        $5 \times 10^{-3}$ && $4.27 \times 10^{-2}$ &&$0.71$  \\
        $2.5 \times 10^{-3}$ &&$2.23 \times 10^{-2}$&&$0.94$ \\
        $1.25 \times 10^{-3}$&&$9.91\times 10^{-3}$&&$1.17$   \\
        $6.25 \times 10^{-4}$ && $2.86 \times 10^{-3}$&&$1.79$  \\
        $3.125 \times 10^{-4}$&& $1.33 \times 10^{-3}$&&$1.10$  \\
        \hline
    \end{tabular}
    \caption{Maximum errors with the convergence rate \\ related to the time discretization.}\label{Table-1}
\end{table}

\begin{figure}[h!]
    \centering
    \includegraphics[width=12cm,height=6cm]{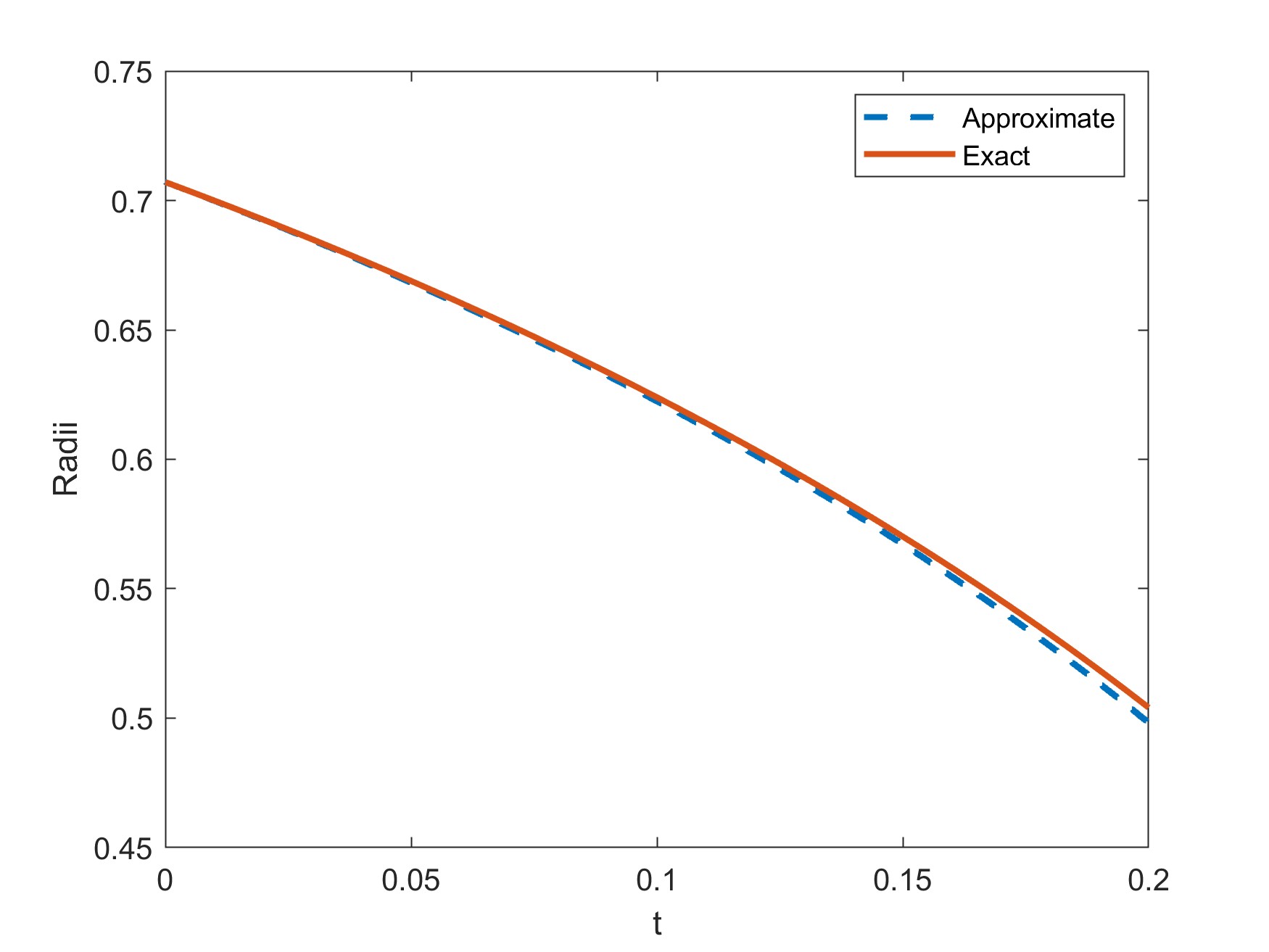}
    \caption{Exact and approximate values of $r(t)$ at different time levels up to $0.2$.}
    \label{Approximate-Exact-Radii}
\end{figure}

\begin{table}[h!]
    \centering
    \begin{tabular}{ccccccccccccccccccc}
        \hline
        \hline
        $N$ && $\ell_{\infty}$ && Orders \\
        \hline
        $1681$ &&$3.78 \times 10^{-2}$&& $-$ \\
        $2601$ && $1.96 \times 10^{-2}$ &&$3.02$  \\
        $3721$ &&$1.34 \times 10^{-2}$&&$2.10$ \\
        $5041$&&$7.69\times 10^{-3}$&&$3.68$   \\
        \hline
    \end{tabular}
    \caption{Maximum errors with the convergence rate \\ related to the spatial discretization.}\label{Table-2}
\end{table}
In addition, Figure \ref{Approximate-Exact-Radii} demonstrates the exact and approximate values of the radii $r$ in time, where $6561$ PTS points and $\Delta t=10^{-4}$ are used.

\begin{figure}[h!]
    \centering
    \includegraphics[width=5cm,height=4.5cm]{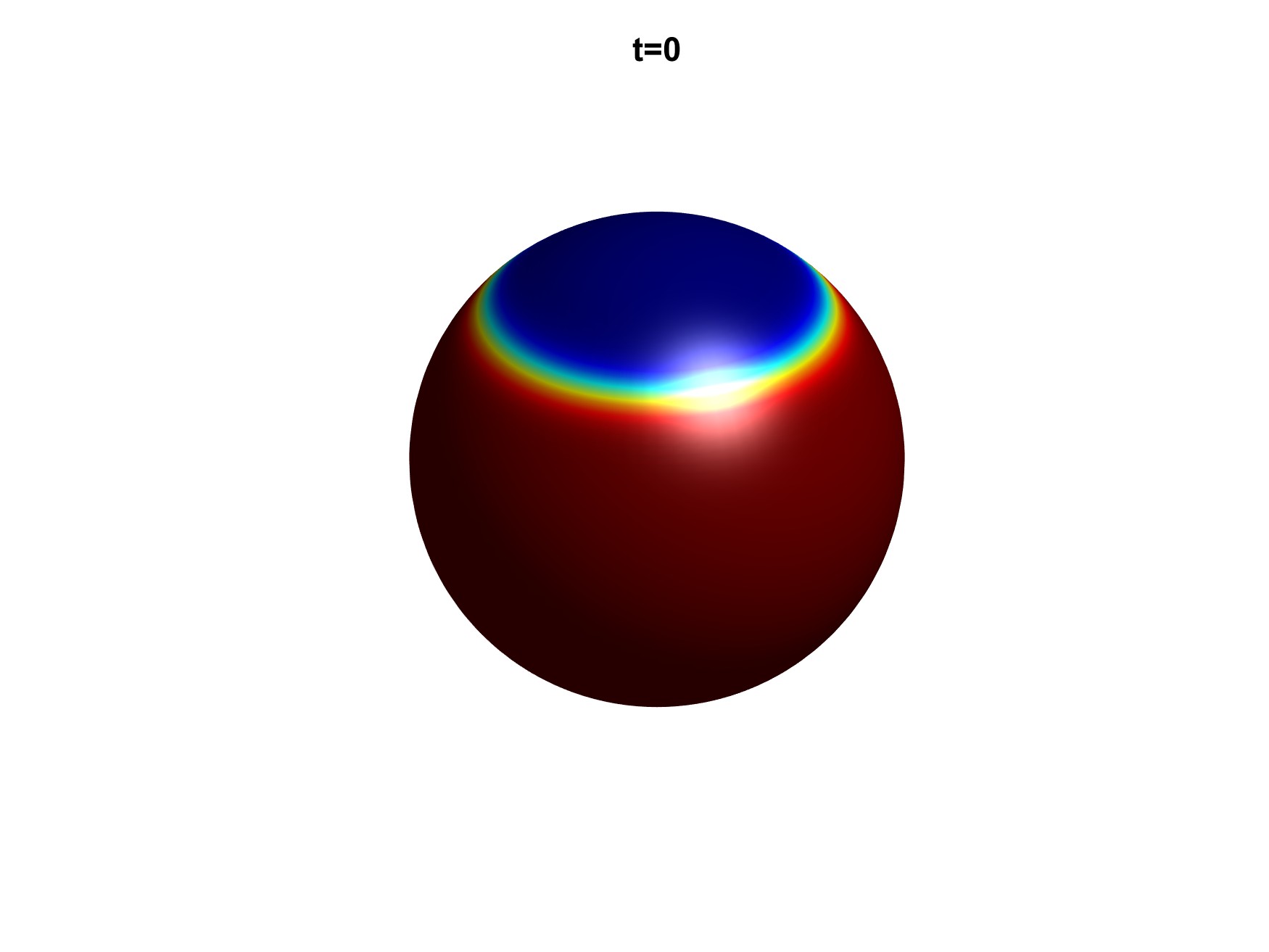}
    \includegraphics[width=5cm,height=4.5cm]{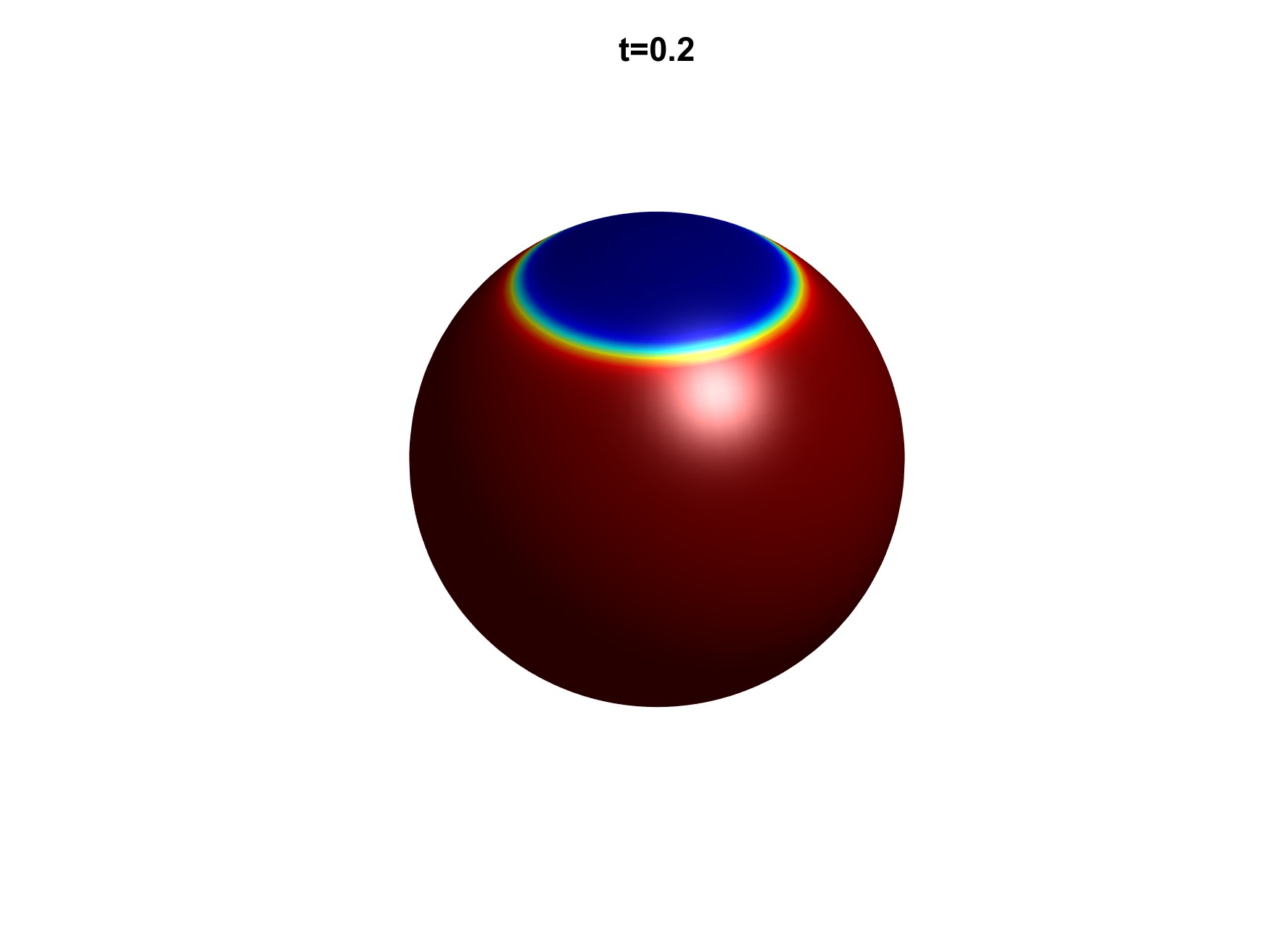}
    \includegraphics[width=5cm,height=4.5cm]{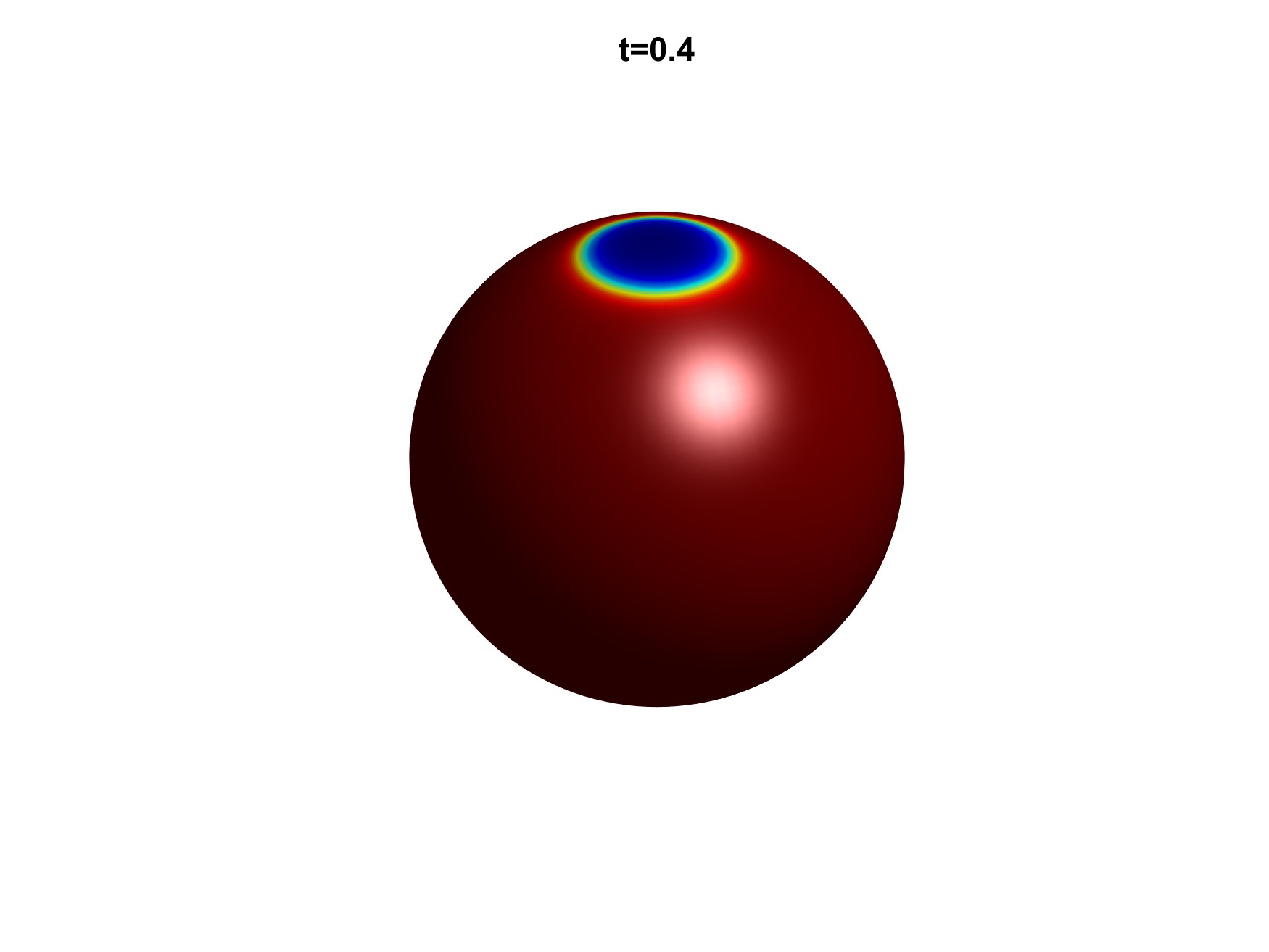}
    \caption{The initial condition with the approximate solution of $u$ at $t=0.2$ and $t=0.4$ via $6561$ PTS points.}
    \label{profile-AC-2}
\end{figure}

\begin{figure}[h!]
    \centering
    \includegraphics[width=5cm,height=4.5cm]{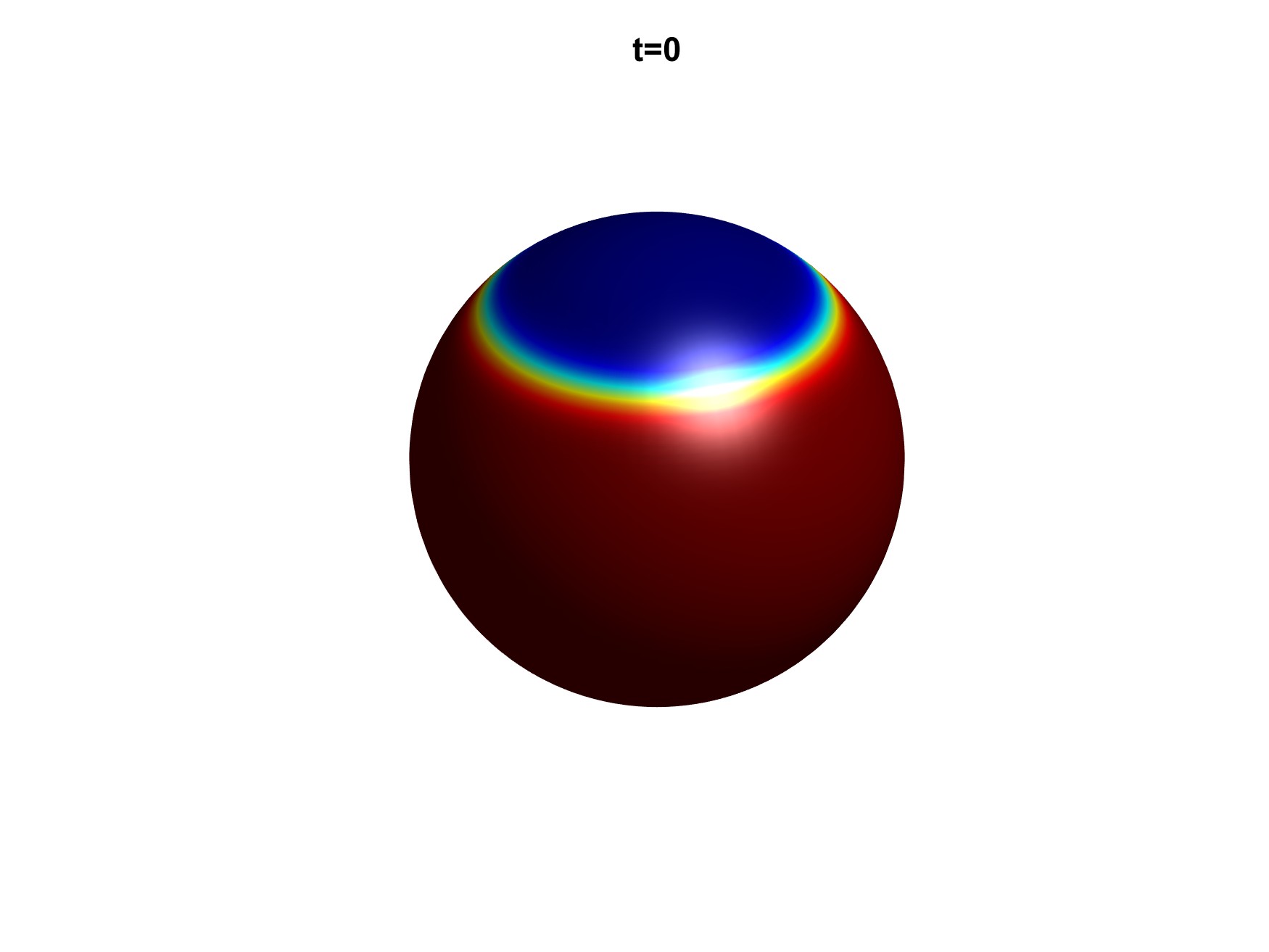}
    \includegraphics[width=5cm,height=4.5cm]{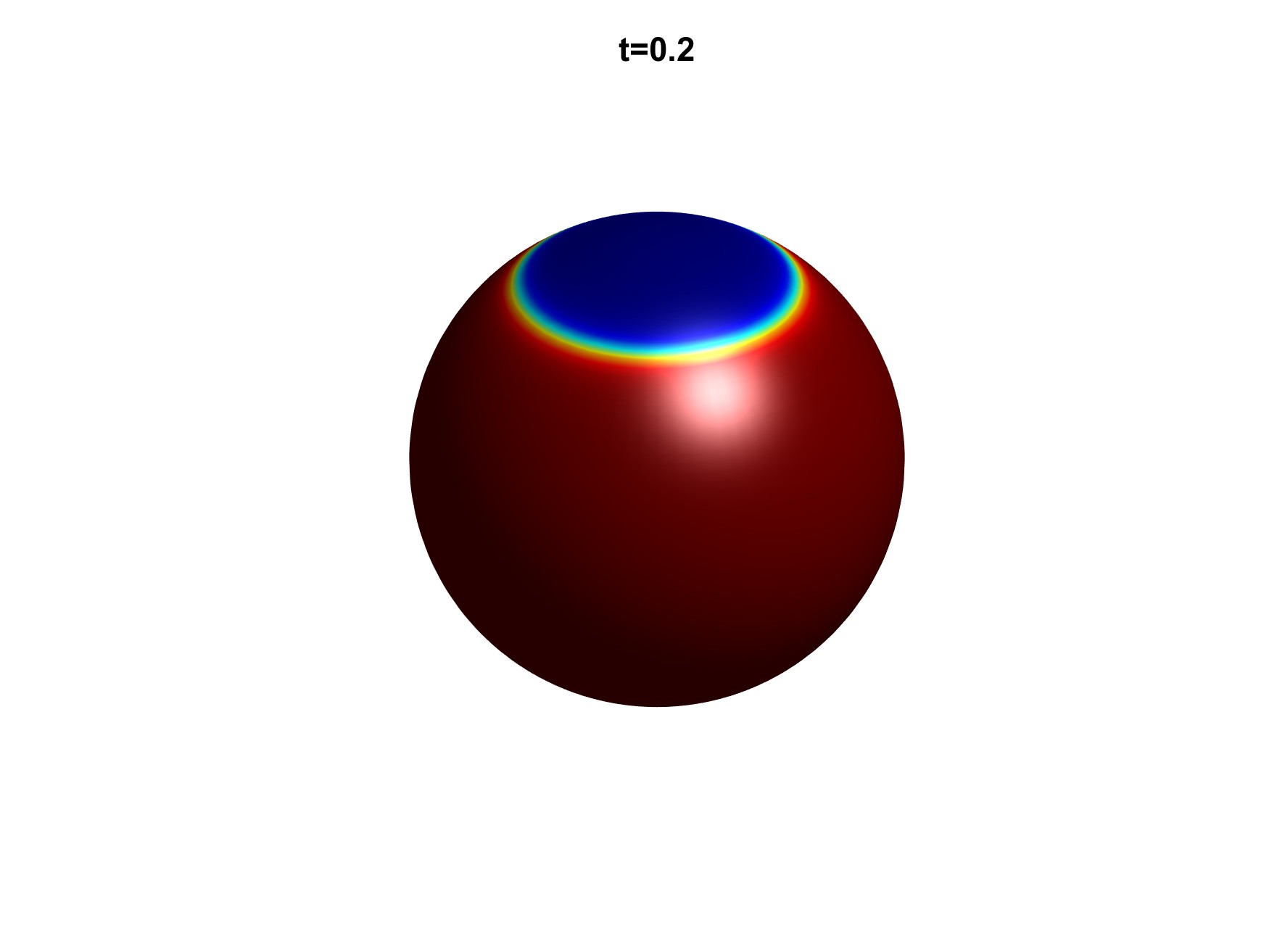}
    \includegraphics[width=5cm,height=4.5cm]{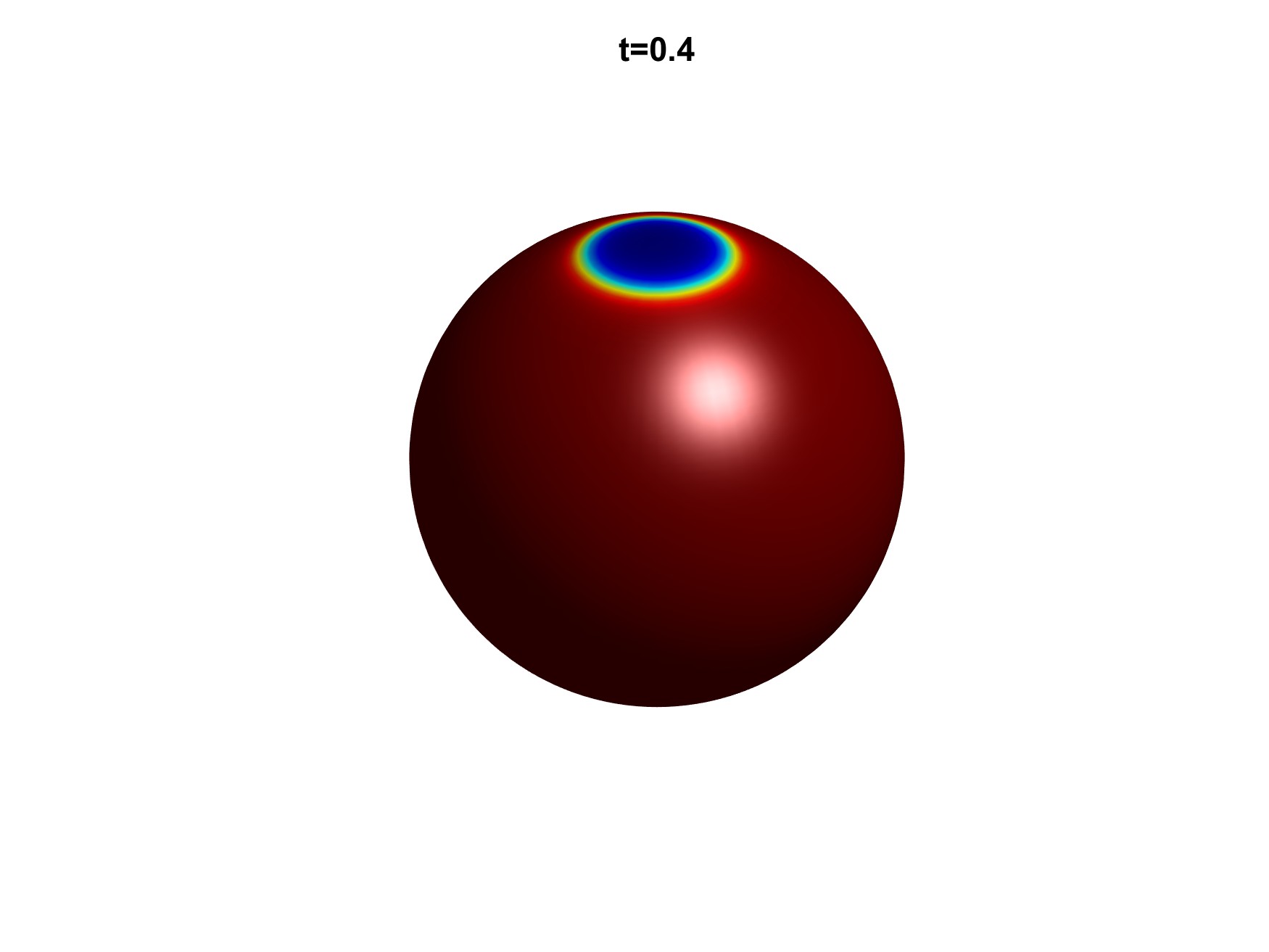}
    \caption{The initial condition with the approximate solution of $u$ at $t=0.2$ and $t=0.4$ via $6561$ ME points.}
    \label{profile-AC-3}
\end{figure}

Besides, to show the temporal evolution of $u$, we have chosen $N$ and $\Delta t$ to be $6561$ and $10^{-2}$, respectively, where both PTS and ME points are utilized. The desired results are drawn in Figures \ref{profile-AC-2} and \ref{profile-AC-3}, which follow from those reported in \cite{Choi,MohDehMir}.

\subsection{Test 2}
In the second test, we present the numerical results related to the CAC equation (\ref{e1-CAC}) on $\mathbb{S}^2$ with $\epsilon=0.0192$ and the initial condition as follows \cite{kim-CAC}:
\newpage
\begin{align}\label{initial-CAC}
u (x,y,z,0) =
\left\{ \begin{array}{ll}
\tanh\Big(\frac{0.8-\sqrt{x^2+y^2}}{\sqrt{2}{\epsilon}}\Big),&z \geq 0,\\
\tanh\Big(\frac{0.5-\sqrt{x^2+y^2}}{\sqrt{2}{\epsilon}}\Big),& z<0,
\end{array} \right.
\end{align}
which shows two spherical caps with radii $0.8$ and $0.5$ located at the north and south poles, respectively. The interested reader is also referred to \cite{kim-CAC} for more information on this test problem.

\begin{figure}[h!]
\centering
\includegraphics[width=7cm,height=5cm]{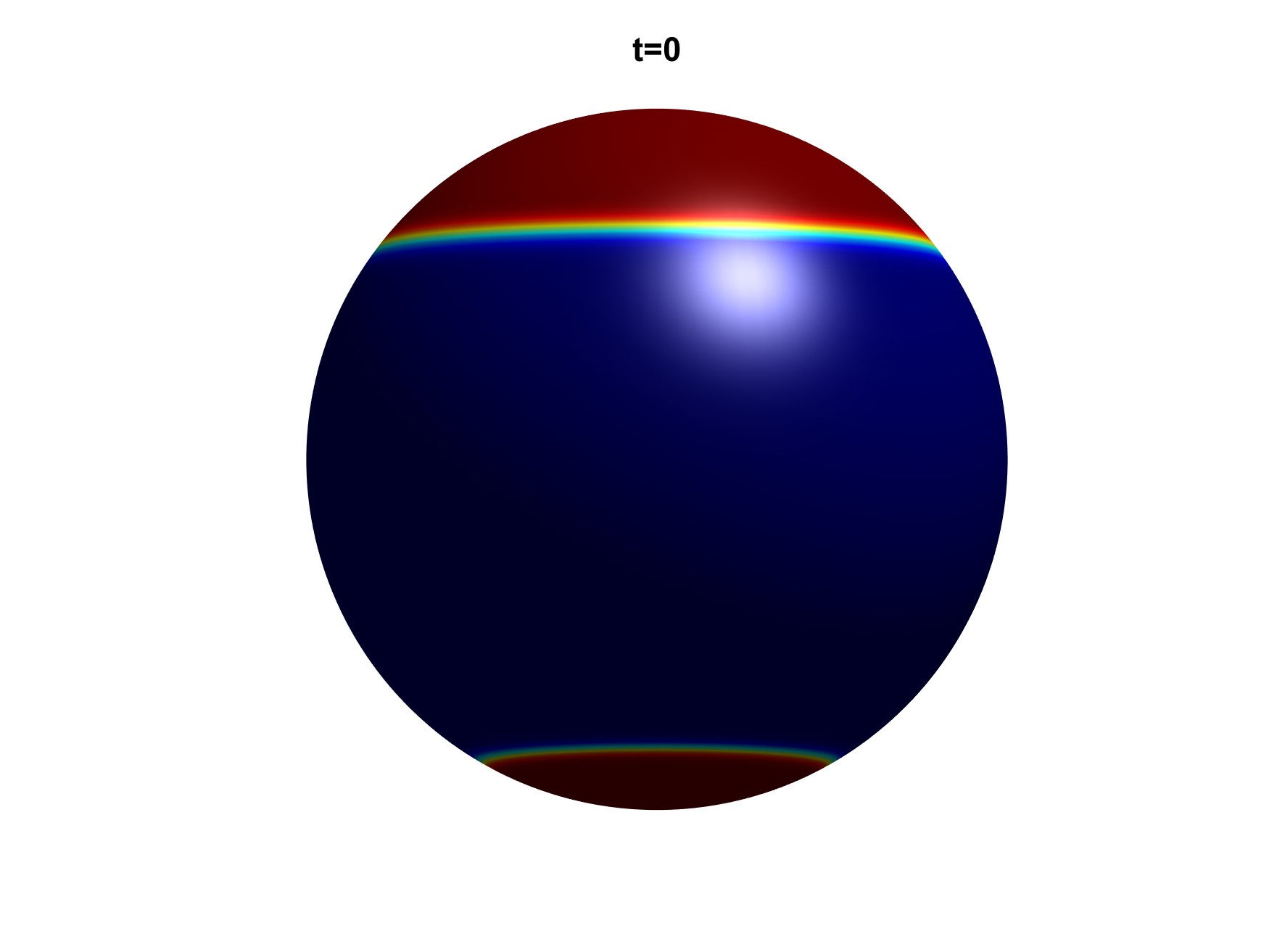}
\includegraphics[width=7cm,height=5cm]{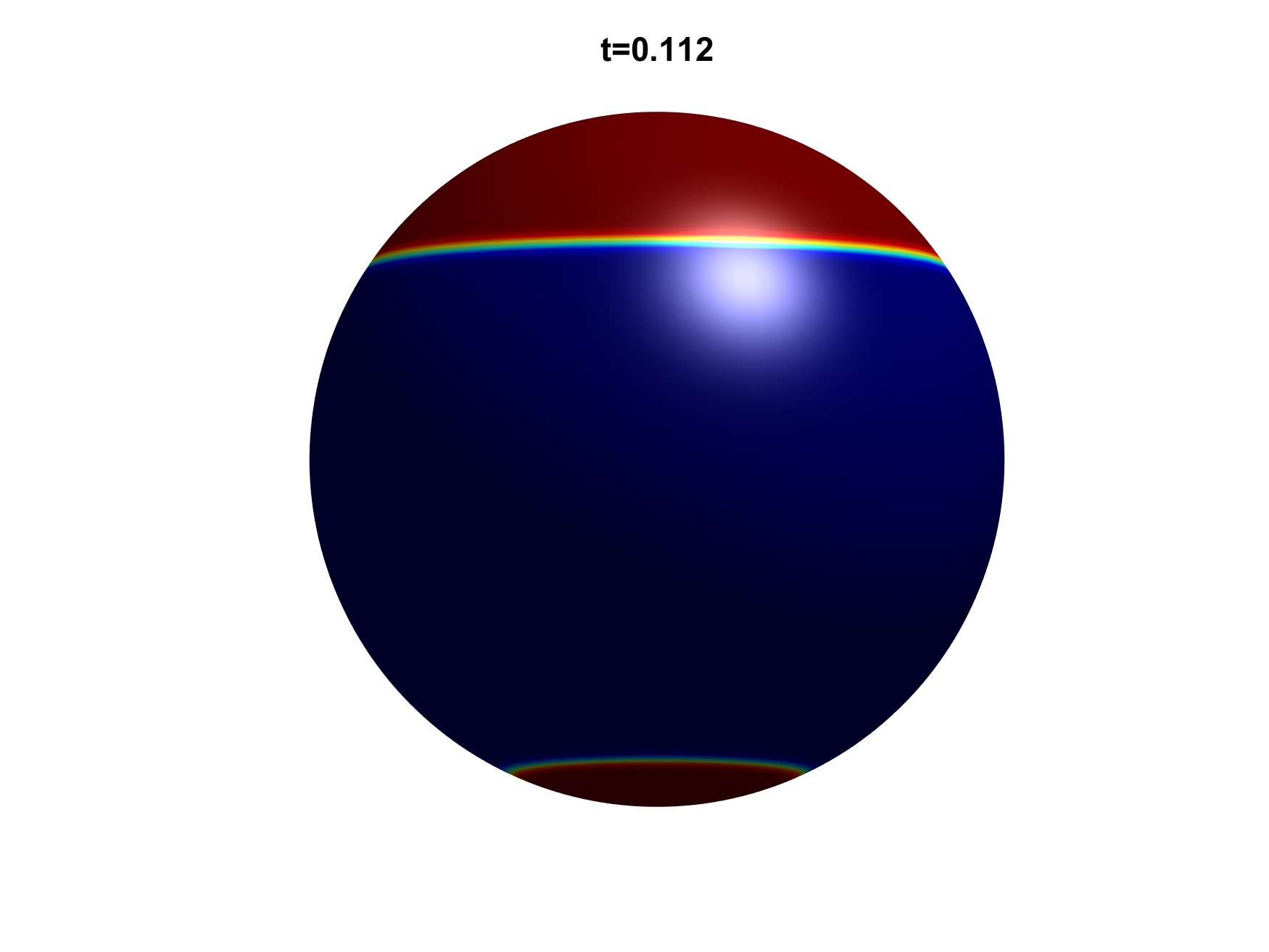}
\includegraphics[width=7cm,height=5cm]{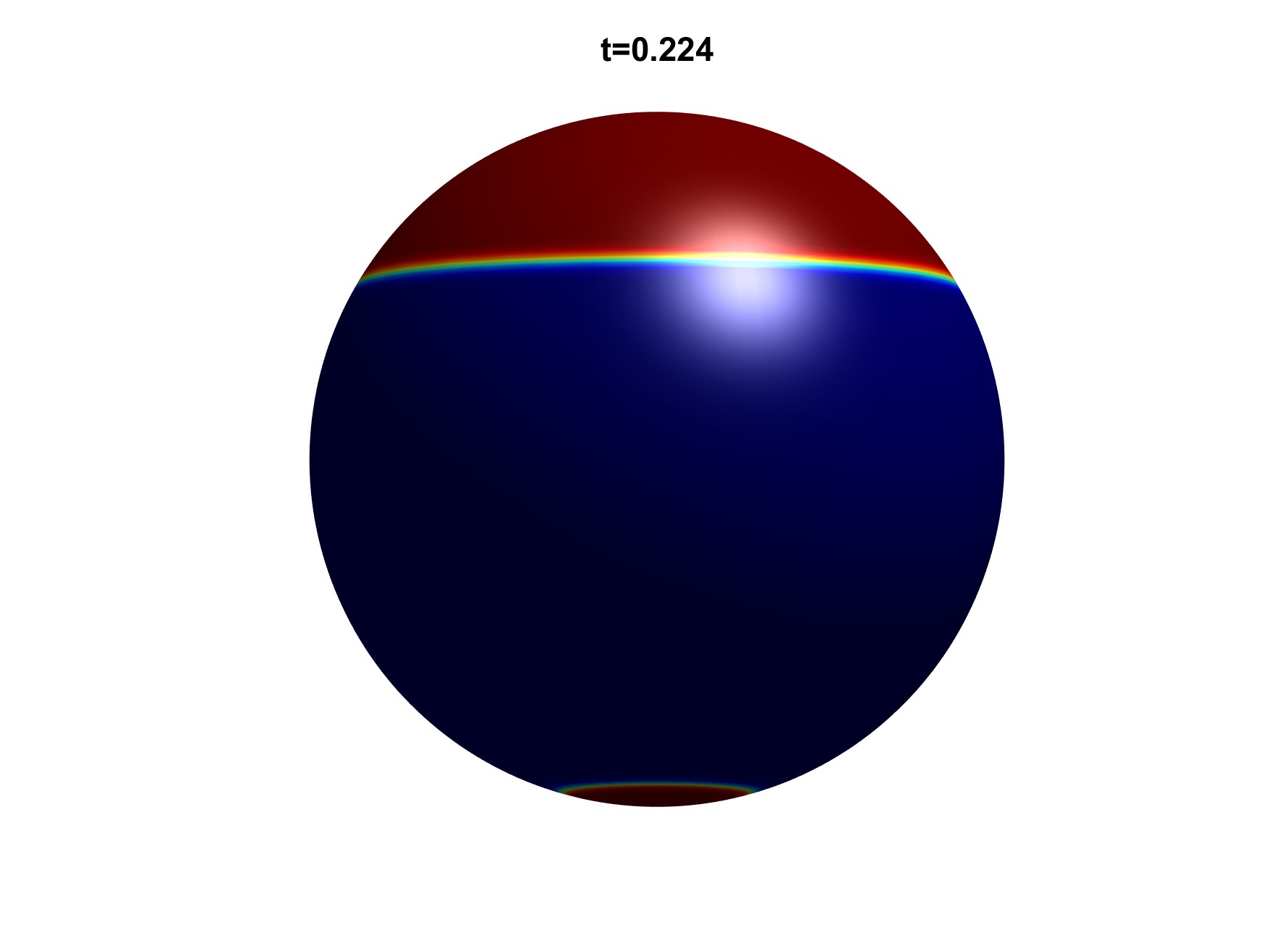}
\includegraphics[width=7cm,height=5cm]{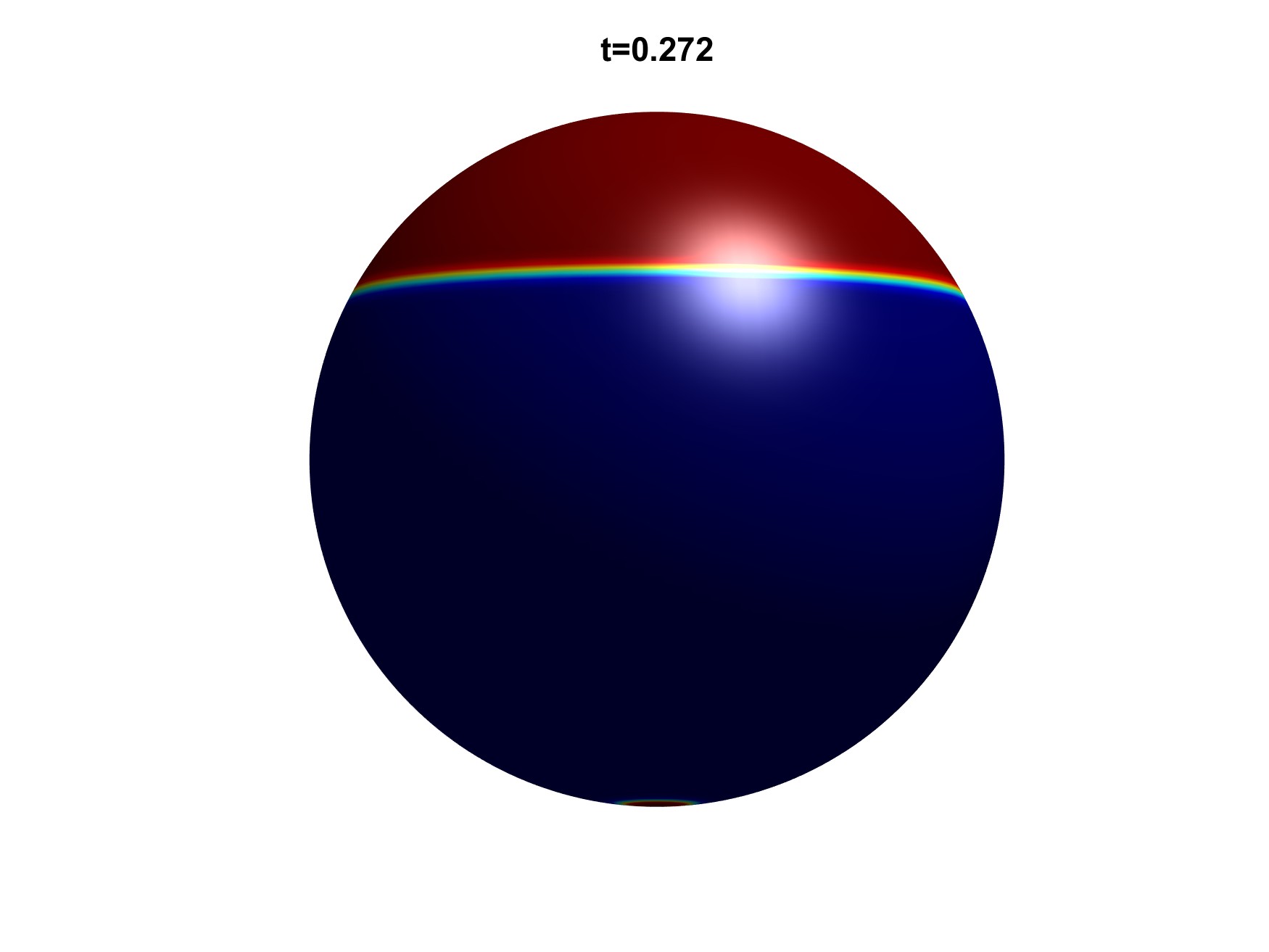}
\caption{The initial condition with the approximate solution of $u$ at $t=0.112$, $t=0.224$ and $t=0.272$ via $25600$ PTS points.}
\label{profile-error-2}
\end{figure}

\begin{figure}[h!]
\centering
\includegraphics[width=14cm,height=8cm]{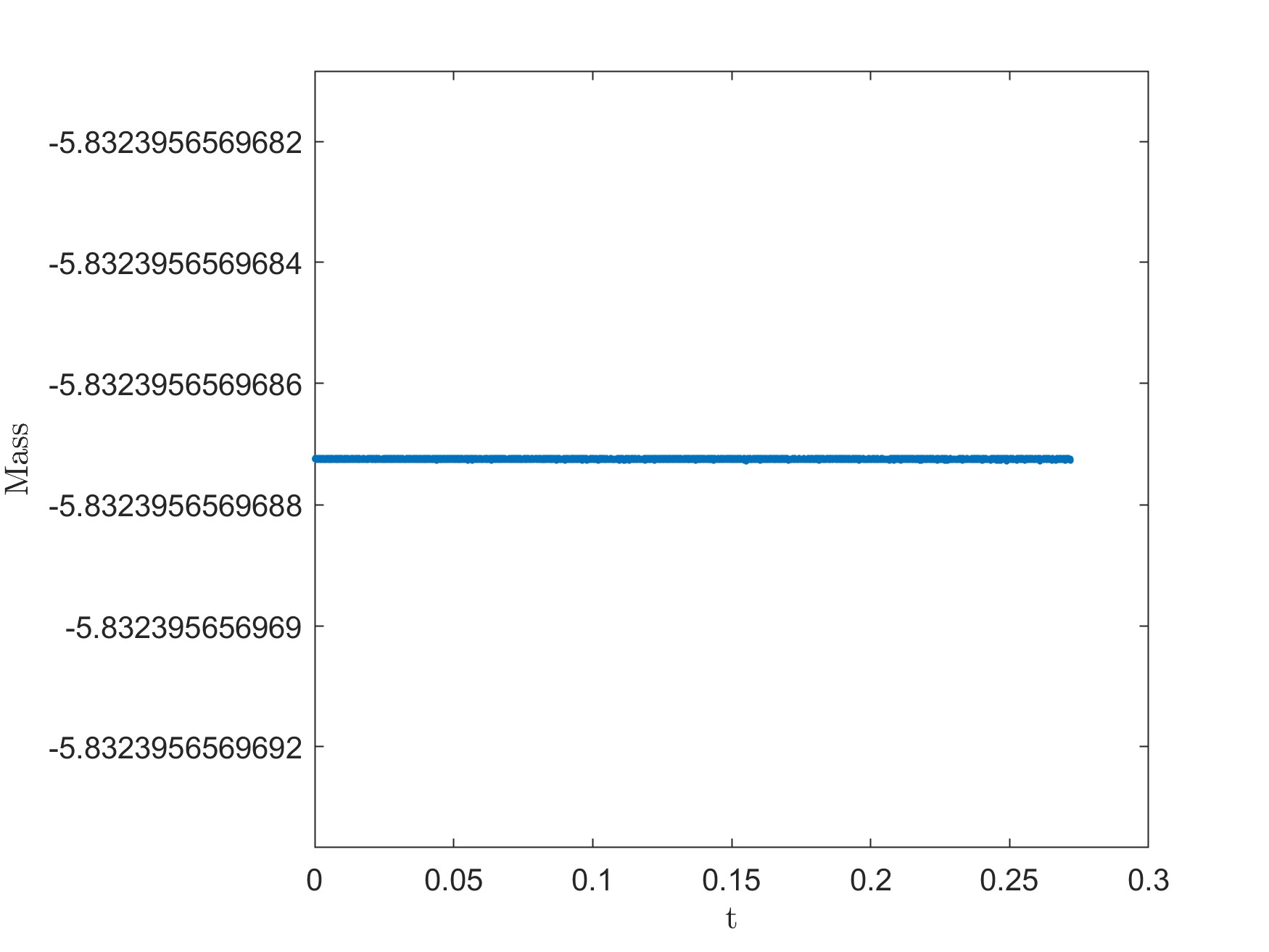}
\caption{Approximate mass values over time.}
\label{profile-Mass-1}
\end{figure}

\begin{figure}[h!]
    \centering
    \includegraphics[width=8cm,height=6cm]{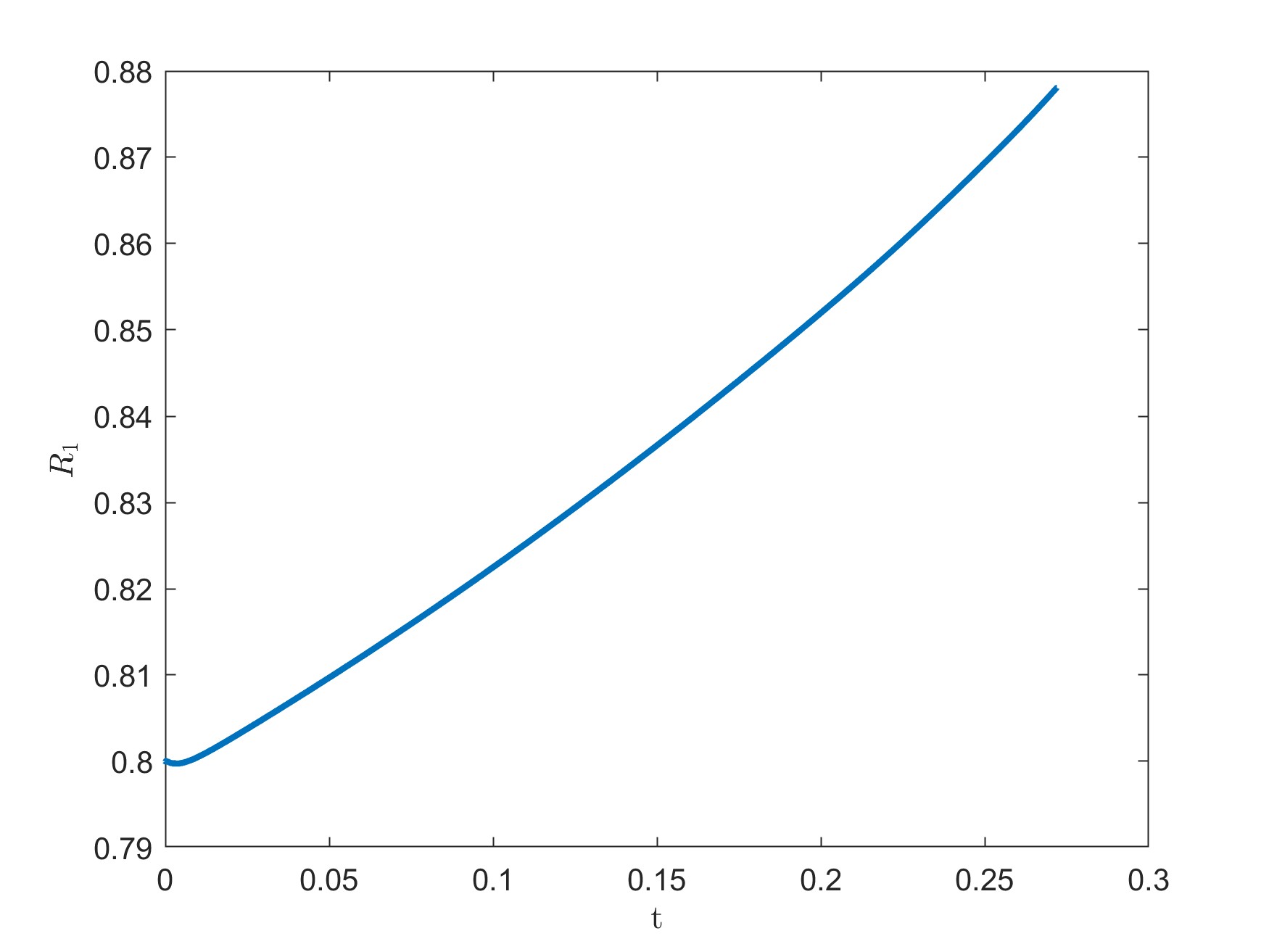}
    \includegraphics[width=8cm,height=6cm]{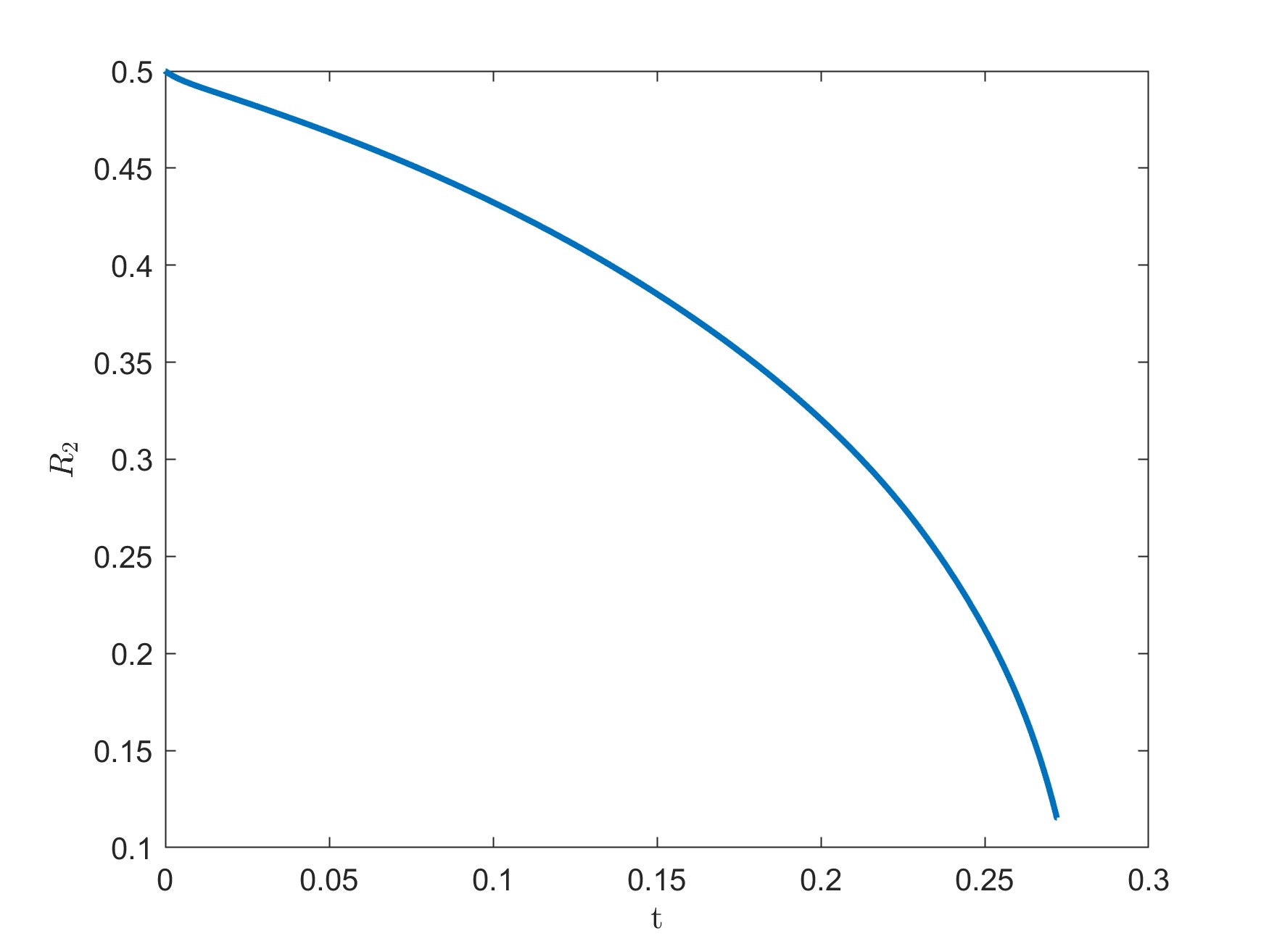}
    \caption{Temporal evolution of the radii $R_{1}$ (left) and $R_{2}$ (right).}
    \label{profile-Radii}
\end{figure}

Here, the PTS points on the unit sphere are used for presenting simulation results. Besides, the radius of each spherical cap is chosen to be $\delta=6m h$, where $m=2$ and $h=1/\sqrt{N}$ with $N=25600$. The time step $\Delta t$ is also fixed, and its value is $4.5 \times 10^{-5}$.
In Figure \ref{profile-error-2}, the initial value with the approximate solution to $u$ is plotted at $t=0.112$, $0.224$, and $0.272$. As can be observed here, the larger spherical cap grows, while its smaller counterpart shrinks until it disappears. However, the same results are obtained in \cite{kim-CAC}, but using more points, i.e., $2803221$, where the quasi-Neumann boundary condition on the narrow band domain using the closest point approach has also been considered. 
Besides, the approximate values of mass in time are displayed in a Figure \ref{profile-Mass-1}. 
In Figure \ref{profile-Radii}, we also report the approximate values of radii $R_{1}$ and $R_{2}$ in time.
\newpage

\section{Concluding remarks}\label{sect-Conclusion}
We have defined a new RKHS whose kernel is constructed from the dual representation of the approximate solution obtained by applying the proposed LS-SVR based on the GMLS technique to the approximation problem on spheres. The connection of the presented learning algorithm to the smoothing problem and KRR has also been investigated in its associated RKHS.  
As an application, two phase-field equations on spheres, i.e., AC and CAC models have been solved via the operator-splitting approach in time and the proposed LS-SVR in space. For doing so, we first discretize both mathematical models using a first-order time-splitting approach, resulting in linear and nonlinear problems. Then, the LS-SVR with the GMLS procedure is applied to the linear problem per time step to deal with the approximate solution. Hence, using the KKT conditions, the unique approximate intermediate solution is computed by solving a linear system of algebraic equations in time, where the coefficient matrix is large, sparse, and symmetric positive definite, and yields the unique numerical solution per time step. Accordingly, the preconditioned conjugate gradient (PCG) algorithm, where an appropriate preconditioner, i.e., ILU(0) factorization of the coefficient matrix is proposed, is executed. After that, the approximate intermediate solution obtained at this stage is directly substituted into the semi-discrete scheme of the nonlinear problem to get the numerical solution. Finally, 
we have presented the results to illustrate the claimed error bounds at points $X$. Moreover, some presented simulation results confirm the suitability of this development in finding stable and accurate numerical solutions for the studied phase-field models on spheres.

\section*{Appendix A}
\begin{lemma}
    If $u(\bm{x},t)$ is the solution of Eq. (\ref{e1-CAC}). Then, it satisfies the following mass conservation property.
    \begin{equation}
    \dfrac{d}{dt}\displaystyle \int_{\mathbb{S}^2}u(\bm{x},t) d\bm{x}=0,\quad t>0.
    \end{equation}
    \begin{proof}
        By multiplying both sides of Eq. (\ref{e1-CAC}) to the test function $1$ in $L_{2}(\mathbb{S}^2)$, we have
        \begin{align}\label{proof-1-Appendix}
         (\dfrac{{\partial u (\bm{x},t)}}{{\partial t}},1)_{L_{2}(\mathbb{S}^2)} =({\Delta _{\mathbb{S}^2}}u
        (\bm{x},t),1)_{L_{2}(\mathbb{S}^2)} - (\dfrac{{F^{\prime}(u (\bm{x},t))}}{{{{\epsilon}
                    ^2}}},1)_{L_{2}(\mathbb{S}^2)}+\beta(t)
                (\sqrt{F(u(\x,t))},1)_{L_{2}(\mathbb{S}^2)}.
        \end{align}
        Since $\mathbb{S}^2$ is a smooth, $2$-dimensional compact manifold, applying the Green-Beltrami (cf. \cite{Atkinson-2}) identity yields
        \begin{align}
        ({\Delta _{\mathbb{S}^2}}u
        (\bm{x},t),1)_{L_{2}(\mathbb{S}^2)}=({\nabla _{\mathbb{S}^2}}u
        (\bm{x},t),{\nabla _{\mathbb{S}^2}} 1)_{\bm{L}_{2}(\mathbb{S}^2)}=0.
        \end{align}
        On the other hand, according to the definition $\beta(t)$, the second and third terms on the right hand side of (\ref{proof-1-Appendix}) can be written as follows:
        \begin{align}
        &- (\dfrac{{F^{\prime}(u (\bm{x},t))}}{{{{\epsilon}
                    ^2}}},1)_{L_{2}(\mathbb{S}^2)}+\beta(t)
        (\sqrt{F(u(\x,t))},1)_{L_{2}(\mathbb{S}^2)}=\\ \nonumber
        &-\displaystyle \int_{\mathbb{S}^2} \dfrac{{F^{\prime}(u (\bm{x},t))}}{{{{\epsilon}
                    ^2}}}d\bm{x}+\Big(
        \dfrac{\displaystyle \int_{\mathbb{S}^2}F^{\prime}(u(\bm{x},t))
            ~d\bm{x}}{\epsilon^2 \displaystyle \int_{\mathbb{S}^2}\sqrt{F(u(\bm{x},t))}~d\bm{x}}\Big)
        \displaystyle\int_{\mathbb{S}^2}\sqrt{F(u(\bm{x},t))}~d\bm{x}=0,
        \end{align}
        which completes the proof.
    \end{proof}

\end{lemma}

\begin{lemma}
    The time discretization (\ref{e5-CAC})-(\ref{e7-CAC}) is mass conservative, i.e.,
    \begin{align}
    (u^{n+1},1)_{L_{2}(\mathbb{S}^2)}=(u^{0},1)_{L_{2}(\mathbb{S}^2)},\quad \forall n \geq 0.
    \end{align}
    \begin{proof}
        According to relation (\ref{e7-CAC}), for $n \geq 0$, we can write
        \begin{align}
        (u^{n+1},1)_{L_{2}(\mathbb{S}^2)}&=
        (u^{n+1^{**}}+\Delta t\beta^{n+1^{**}} \sqrt{F(u^{n+1^{**}})},1)_{L_{2}(\mathbb{S}^2)}\\ \nonumber
        &=(u^{n+1^{**}},1)_{L_{2}(\mathbb{S}^2)}+\Delta t \beta^{n+1^{**}}(\sqrt{F(u^{n+1^{**}})},1)_{L_{2}(\mathbb{S}^2)}\\ \nonumber
        &=(u^{n+1^{**}},1)_{L_{2}(\mathbb{S}^2)}+\Delta t
        (\dfrac{(u^{0}-u^{n+1^{**}},1)_{L_{2}(\mathbb{S}^2)}}{\Delta t (\sqrt{F(u^{n+1^{**}})},1)_{L_{2}(\mathbb{S}^2)}})(\sqrt{F(u^{n+1^{**}})},1)_{L_{2}(\mathbb{S}^2)}\\ \nonumber
        &=(u^{0},1)_{L_{2}(\mathbb{S}^2)},
        \end{align}
    where in the third equality, we have used relation (\ref{e8-CAC}).
    \end{proof}
\end{lemma}

\section*{Appendix B}
Here are the developed numerical schemes applied to Eqs. \eqref{e1}
and \eqref{e1-CAC} for the reader's convenience, where Algorithm
\ref{Algorithm1} refers to the GMLS approximation. Algorithms
\ref{Algorithm2} and \ref{Algorithm3} propose the detailed
procedures.

\begin{algorithm}
\caption{The computational algorithm of the GMLM approach on
$\mathbb{S}^2$ for approximating the Laplace-Beltrami
operator.}\label{Algorithm1}
 \begin{algorithmic}
\State \textbf{Input:} $X=[\bm{x}_{1},...,\bm{x}_{N}]^{T} \in
\mathbb{R}^{N \times 3}$ as trial points on $\mathbb{S}^2$,
$Y=[\bm{y}_{1},...,\bm{y}_{M}]^{T} \in \mathbb{R}^{M \times 3}$ as
trial points on $\mathbb{S}^2$, $\delta$ as a radius of the
spherical
 cap;
 \State \textbf{Output:} $\bm{a}(\bm{x})$ and $
\texttt{L}\bm{a}(\bm{x})$ at test points $Y$ to construct sparse matrices $A_{X}$ and $L_{X}$;
\State Set the matrices $A_{X}$ and $L_{X}$ to be zeros of size $M\times N$;
\For {$i=1,...,M$}
 \State Construct $I(\bm{y}_{i})$ with respect to $X$;
 \State
Compute the matrix $W$ by choosing the weight function
(\ref{Weight-Function});
 \State Compute the local matrix $P^{T}W(\bm{y}_{i})P$ with respect to trial points in
 $I(\bm{y}_{i})$;
 \State Compute the vectors $\bm{Y}^{T}(\bm{y}_{i})$ and $\texttt{L}\bm{Y}^{T}(\bm{y}_{i})$ (see relation (\ref{GMLS-Sphere2}));
 \State Compute the vectors
$\bm{a}(\bm{y}_{i})$ and $\texttt{L}\bm{a}(\bm{y}_{i})$
  due to Eqs. (\ref{GMLS-Sphere-new-1}) and (\ref{GMLS-Sphere6}) and using the QR factorization of $\sqrt{W}P$
  to find $C=:WP\Big(P^{T}W(\bm{y}_{i})P\Big)^{-1}$;
  \State $A_{X}(i,I(\bm{y}_{i}))\leftarrow \bm{a}(\bm{y}_{i})$;
  \State $L_{X}(i,I(\bm{y}_{i}))\leftarrow \texttt{L}\bm{a}(\bm{y}_{i})$;
 \EndFor
 \State
\end{algorithmic}
\end{algorithm}

\begin{algorithm}
    \caption{The computational algorithm of the numerical scheme applied
    to
    the AC equation (\ref{e1}) on $\mathbb{S}^2$.}\label{Algorithm2}
    \begin{algorithmic}[1]
    \State \textbf{Input:}$X=[\bm{x}_{1},...,\bm{x}_{N}]^{T}
                \in \mathbb{R}^{N \times 3}$ as data
                points on $\mathbb{S}^2$, $\delta>0$ as the spherical cap,
                $\epsilon$, $\gamma$, $T$ as final time, and the time step
                $\Delta t$;
                \State \textbf{Output:} Approximate solution $U_{X}$ at the final time $T$;
                \State Construct the sparse matrices $A_{X}$ and $A_{L}$ of size $N \times N$ due to Algorithm \ref{Algorithm1};
                \State Form the matrix $B=A_{L}A^{T}_{L}+\frac{1}{\gamma}I_{N}$;
                \State Apply the RCM algorithm to the coefficient matrix $B$, called $B_{R}$;
                \State Use ILU(0) factorization to the matrix $B_{R}$
                to construct the proper preconditioner;
                \State $t\leftarrow 0$;
                \State $n\leftarrow 0$;
                \State Set the vector of initial condition at points $X$, i.e.,
                 ${U}_{X}^{0} \in \mathbb{R}^{N \times 1}$;
                \While {$t \leq T$}
                \State Find $\bm{\alpha}^{*}_{n}$ by applying the RCM algorithm to the linear system
                (\ref{LS-SVM-7}), and solve it via the PCG algorithm;
                \State Compute $\bm{w}^{*}_{n}$ from the first relation in
                (\ref{LS-SVM-7-1});
                \State Compute the intermediate solution $U^{n*}_X$
                via (\ref{LS-SVM-8});
                \State Find the approximate solution $U^{n+1}_X$
                using (\ref{LS-SVM-9});
                       \State $n\leftarrow n+1$;
                \State $t \leftarrow n \Delta t$;
                \EndWhile
                \State ${U}_{X} \leftarrow {U}^{K}_{X}$;
    \end{algorithmic}
\end{algorithm}

\begin{algorithm}
    \caption{The computational algorithm of the numerical scheme applied
    to
    the CAC equation (\ref{e1-CAC}) on $\mathbb{S}^2$.}\label{Algorithm3}
    \begin{algorithmic}[1]
    \State \textbf{Input:}$X=[\bm{x}_{1},...,\bm{x}_{N}]^{T}
                \in \mathbb{R}^{N \times 3}$ as data
                points on $\mathbb{S}^2$, $\delta>0$ as the spherical cap,
                $\epsilon$, $\gamma$, $T$ as final time, and the time step
                $\Delta t$;
                \State \textbf{Output:} Approximate solution $U_{X}$ at the final time $T$;
                \State Construct the sparse matrices $A_{X}$ and $A_{L}$ of size $N \times N$ due to Algorithm \ref{Algorithm1};
                \State Form the matrix $B=A_{L}A^{T}_{L}+\frac{1}{\gamma}I_{N}$;
                \State Apply the RCM algorithm to the coefficient matrix $B$, called $B_{R}$;
                \State Use ILU(0) factorization to the matrix $B_{R}$
                to construct the proper preconditioner;
                \State $t\leftarrow 0$;
                \State $n\leftarrow 0$;
                \State Set the vector of initial condition at points $X$, i.e.,
                 ${U}_{X}^{0} \in \mathbb{R}^{N \times 1}$;
                \While {$t \leq T$}
                \State Find $\bm{\alpha}^{*}_{n}$ by applying the RCM algorithm to the linear system
                (\ref{LS-SVM-7-CAC}), and solve it via the PCG algorithm;
                \State Compute $\bm{w}^{*}_{n}$ from the first relation in
                (\ref{LS-SVM-7-1-CAC});
                \State Compute the intermediate solution $U^{n^{*}}_X$
                via (\ref{LS-SVM-8-CAC});
                \State Find the approximate solution $U^{{n+1}^{**}}_X$
                using (\ref{LS-SVM-9-CAC});
                \State Compute the coefficient $\beta^{n+1^{**}}$ via
                the numerical integration (\ref{Integration});
                \State Find the approximate solution $U^{n+1}_X$
                via (\ref{LS-SVM-10-CAC});
                \State $n\leftarrow n+1$;
                \State $t \leftarrow n \Delta t$;
                \EndWhile
                \State ${U}_{X} \leftarrow {U^{K}_X}$;
    \end{algorithmic}
\end{algorithm}

\end{document}